\documentclass[10pt]{IEEEtran}

\usepackage{cite}
\usepackage{url}
\usepackage{array}
\usepackage{amsfonts}
\usepackage[cmex10]{amsmath}
\usepackage{amssymb}

\makeatletter
\@addtoreset{equation}{section}
\makeatother
\newtheorem{theorem}{Theorem}[section]
\newtheorem{lemma}[theorem]{Lemma}
\newtheorem{corollary}[theorem]{Corollary}
\newtheorem{proposition}[theorem]{Proposition}
\newtheorem{fact}[theorem]{Fact}
\newtheorem{definition}[theorem]{Definition}

\newtheorem{remark}[theorem]{Remark}

\renewcommand{\IEEEQED}{\IEEEQEDopen}

\begin{document}

\title{Linear Nonbinary Covering Codes and Saturating Sets in
Projective Spaces}
\author{Alexander~A.~Davydov, Massimo~Giulietti, Stefano~Marcugini,
and~Fernanda~Pambianco
\thanks{Manuscript received xxxx 00, 2009; revised xxxx 00, xxxx.
The material in this
paper was presented in part at the XI International Workshop on
Algebraic and Combinatorial Coding Theory, ACCT2008, Pamporovo,
Bulgaria, June 2008, and at International Workshop ``Coding
Theory Days in St. Petersburg'', St. Petersburg, Russia, Oct.
2008.}
\thanks{A. A. Davydov is with the Institute for Information Transmission
Problems, Russian Academy of Sciences, GSP-4, Moscow, 127994,
Russian Federation (e-mail: adav@iitp.ru).}
\thanks{M.\ Giulietti, S.\ Marcugini, and F.\ Pambianco are with the Department
of Mathematics and Informatics, Perugia University, Perugia,
06123, Italy (e-mail: giuliet@dipmat.unipg.it;
gino@dipmat.unipg.it; fernanda@dipmat.unipg.it).}
\thanks{Communicated by xxxxxxxxxxxxxxxxxxx.}
\thanks{Digital Object Identifier 00.0000/xxx.0000.000000}}
\markboth{IEEE Transactions on Information Theory,~vol.~00,
no.~0,~xxxxxxxx~0000}{}
\IEEEpubid{0000--0000/00\$00.00~\copyright~0000 IEEE}
\maketitle

\begin{abstract}
Let $\mathcal{A}_{R,q}$ denote a family of covering codes, in
which the covering radius $R$ and the size $q$ of the
underlying Galois field are fixed, while the code length tends
to infinity. The construction of families with small asymptotic
covering densities is a classical problem in the area Covering
Codes.

In this paper, infinite sets of families $\mathcal{A}_{R,q}$,
where $R$ is fixed but $q$ ranges over an infinite set of prime
powers are considered, and the dependence on $q$ of the
asymptotic covering densities of $\mathcal{A}_{R,q}$ is
investigated. It turns out that for the upper limit $\mu
_{q}^{\ast}(R,\mathcal{A}_{R,q})$ of the covering density of
$\mathcal{A}_{R,q}$, the best possibility is
\begin{equation}
\mu
_{q}^{\ast }(R,\mathcal{A}_{R,q})= O(q).\label{ab}
\end{equation}
The main achievement of the present paper is the construction
of {\em optimal} infinite sets of families $\mathcal{A}_{R,q}$,
that is, sets of families such that (\ref{ab}) holds, for any
covering radius $R\ge 2$.

We first showed that for a given $R$, to obtain optimal
infinite sets of families it is enough to construct $R$
infinite families
$\mathcal{A}_{R,q}^{(0)},\mathcal{A}_{R,q}^{(1)},\ldots
,\mathcal{A}_{R,q}^{(R-1)}$ such that, for all $u\geq u_{0}$,
the family $\mathcal{A}_{R,q}^{(\gamma )}$ contains codes of
codimension $r_{u}=Ru+\gamma$ and length $f_{q}^{(\gamma
)}(r_{u})$ where $f_{q}^{(\gamma )}(r)=O(q^{\frac{r-R}{R}})$
and $u_{0}$ is a constant. Then, we were able to construct the
needed families $\mathcal{A}_{R,q}^{(\gamma )}$ for any
covering radius $R\geq 2$, with $q$ ranging over the (infinite)
set of $R$-th powers. A result of independent interest is that
in each of these families $\mathcal{A}_{R,q}^{(\gamma )}$, the
lower limit of the covering density is bounded from above by a
constant independent of $q$.

The key tool in our investigation is the design of new small
saturating sets in projective spaces over finite fields, which
are used as the starting point for the $q^{m}$-concatenating
constructions of covering codes. A new concept of $N$-fold
strong\ blocking set is introduced. As a result of our
investigation, many new asymptotic and finite upper bounds on
the length function of covering codes and on the smallest sizes
of saturating sets, are also obtained. Updated tables for these
upper bounds are provided. An analysis and a survey of the
known results are presented.
\end{abstract}


\renewcommand{\theequation}{\arabic{section}.\arabic{equation}}

\section{Introduction}

Let $F_{q}$ be the Galois field with $q$ elements. Let
$F_{q}^{n}$ be the $n$ -dimensional vector space over $F_{q}.$
Denote by $[n,n-r]_{q}$ a $q$-ary linear code of length $n$ and
codimension (redundancy) $r,$ that is, a subspace of
$F_{q}^{n}$ of dimension $n-r.$ For an introduction to coding
theory, see \cite{M-WSl},\cite{Handbook-codes}.

\IEEEpubidadjcol

The Hamming distance $d(v,c)$ of vectors $v$ and $c$ in
$F_{q}^{n}$ is the number of positions in which $v$ and $c$
differ. The smallest Hamming distance between distinct code
vectors is called the minimum distance of the code. An
$[n,n-r]_{q}$ code with minimum distance $d$ is denoted as an $
[n,n-r,d]_{q}$ code. The sphere of radius $R$ with center $c$
in $F_{q}^{n}$ is the set $\{v:v\in F_{q}^{n},$ $d(v,c)\leq
R\}$.

\begin{definition}
\label{Def1_CoverRad}\emph{i)} The \emph{covering radius }of an $[n,n-r]_{q}$
code is the least integer $R$ such that the space $F_{q}^{n}$ is covered by
spheres of radius $R$ centered on codewords.

\emph{ii)} A linear $[n,n-r]_{q}$ code has \emph{covering radius} $R$ if
every column of $F_{q}^{r}$ is equal to a linear combination of $R$ columns
of a parity check matrix of the code, and $R$ is the smallest value with
such property.
\end{definition}

Definition \ref{Def1_CoverRad}i makes sense for both linear and
nonlinear codes. For linear codes Definitions
\ref{Def1_CoverRad}i and \ref {Def1_CoverRad}ii are equivalent.
An $[n,n-r]_{q}R$ code ($[n,n-r,d]_{q}R$ code, resp.) is an
$[n,n-r]_{q}$ code ($[n,n-r,d]$ code, resp.) with covering
radius~$R$. For an introduction to coverings of vector Hamming
spaces over finite fields, see
\cite{Coh,Handbook-coverings,CohSurvey1985,GrSl}.

The covering problem for codes is that of finding codes with
small covering radius with respect to their lengths and
dimensions. Codes investigated from the point of view of the
covering problem are usually called \emph{covering codes }(in
contrast to error-correcting codes) \cite{GrSl}.

Problems connected with covering codes are considered in
numerous works, see e.g. \cite{KabPan} -- \cite{DavOprep} and
the references therein, the references in
\cite{Coh,Handbook-coverings,CohSurvey1985,GrSl}, and the
online bibliography of \cite {Lobstein}. In this work, we
mainly give references to researches on nonbinary codes; some
papers on binary codes are also mentioned as they contain
useful general ideas. It should be noted that the monographs
\cite {Coh},\cite{Handbook-coverings} mostly deal with binary
covering codes, and that no surveys of nonbinary covering codes
have been recently published. In this work we try to make up
for this deficiency for linear codes; in particular, for
infinite linear code families. We obtain a number of new
asymptotic optimal results, essentially improving the known
estimates for both finite and infinitely growing code lengths.
The description of new results is provided, along with a survey
of the known ones and their updates.

Studying covering codes is a classical combinatorial task. Covering codes
are connected with many areas of information theory and combinatorics, see,
e.g., \cite[Sec.\thinspace 1.2]{Coh} where problems of data compression,
decoding errors and erasures, football pools, write-once memories,
Berlekamp-Gale game, and Caley graphs are mentioned. Covering codes can also
be used in steganography, see \cite{GalKaba},\cite{BierbStegan} and the
references therein. Codes of covering radius 2 and codimension 3 are
relevant, for example, for defining sets of block designs \cite{Bo-Sz-Ti}
and the degree/diameter problem in graph theory \cite{KKKPS}. Covering codes
can be used in databases \cite{PartSumQuer}. There are connections between
covering codes and a popular game puzzle, called \textquotedblleft
Hats-on-a-line\textquotedblright\ \cite{Hats-on-line}. Covering codes can be
also used to construct identifying codes \cite{Identif}.

It should be particularly emphasized that linear covering codes
are deeply connected with \emph{saturating sets} in
\emph{projective spaces} over finite fields. Let $PG(v,q)$ be
the $v$-dimensional projective space over $ F_{q}$. For an
introduction to such spaces and the geometrical objects
therein, see \cite{Hirs,Hirs1,HS1,HS,Sz-survey89,Lang2000DM}.

A set of points $S\subseteq PG(v,q)$ is said to be $\varrho
$\emph{ -saturating} if for any point $x\in PG(v,q)$ there
exist $\varrho +1$ points in $S$ generating a subspace of
$PG(v,q)$ containing $x$, and $\varrho $ is the smallest value
with such property, cf.
\cite{Bartocci,Ughi,DavO4,MPAustr,DMP-JCTA},
\cite{Jan},\cite[Def.\thinspace
1.1]{Dav95},\cite{DMPIEEE},\cite {DFMP-IEEE-LO}.

As usual, by an $n$-set of $PG(v,q)$ we mean a point set of
size $n.$ Homogenous coordinates of points of an
$(R-1)$-saturating $n$-set $K$ in the projective space
$PG(r-1,q)$ can be treated as columns of a parity-check matrix
of an $[n,n-r]_{q}R$ \emph{related }covering code
$\mathcal{C}_{K}$
\cite{Jan},\cite{Dav95},\cite{DavO5},\cite{DMPIEEE},\cite{DavO4},\cite
{DMP-JCTA}.

In the literature, saturating sets are also called
\textquotedblleft saturated sets\textquotedblright\
\cite{Ughi},\cite{Jan},\cite{Dav95},\cite
{DavNBCR2},\cite{kovS92}, \textquotedblleft spanning
sets\textquotedblright\ \cite{BrPlWi}, and \textquotedblleft
dense sets\textquotedblright
\cite{Bo-Sz-Ti},\cite{Bartocci},\cite{GiuTor-Ars04},\cite{Giul-plane}.

Let $V_{q}(n,R)$ be the cardinality of the sphere of radius $R$
in the vector space $F_{q}^{n}.$
\begin{equation}
V_{q}(n,R)=\sum_{i=0}^{R}(q-1)^{i}\binom{n}{i}.  \label{form1_sphere}
\end{equation}
The covering quality of an $[n,n-r(\mathcal{C})]_{q}R$ code $\mathcal{C}$ of
codimension $r(\mathcal{C})$ can be measured by\emph{\ }its \emph{covering
density}
\begin{equation}
\mu _{q}(n,R,\mathcal{C})=\frac{V_{q}(n,R)}{q^{r(\mathcal{C})}}.
\label{form1_covdensity}
\end{equation}
We will write $\mu _{q}(n,R)$ for $\mu _{q}(n,R,\mathcal{C})$
when the code $ \mathcal{C}$ is clear from the context. Note
that $\mu _{q}(n,R,\mathcal{C} )\geq 1$ and equality holds when
$\mathcal{C}$ is a perfect code. From the point of view of
covering problem, the best codes are those with small covering
density.

For fixed parameters $r,R,$ and $q,$ the smaller is the length
$n$ of an $ [n,n-r]_{q}R$ code, the smaller is its covering
density. The\emph{\ length function} $\ell _{q}(r,R)$ is the
smallest length of a $q$-ary linear code with codimension $r$
and covering radius $R$ \cite{BrPlWi},~\cite
{Handbook-coverings}. The \emph{smallest known} length of such
code is denoted by $\overline{\ell }_{q}(r,R).$ Clearly, $\ell
_{q}(r,R)\leq \overline{\ell }_{q}(r,R)$ holds, and the
existence of an $[n,n-r]_{q}R$ code or an $(R-1)$-saturating
$n$-set in $PG(r-1,q)$ implies the upper bound $\overline{\ell
}_{q}(r,R)\leq n.$

\begin{fact}
\label{Fact1_lengthening}If there is an $[n,n-r]_{q}R$ code then there is an
$[n+1,n+1-r]_{q}R$ code.
\end{fact}

One can obtain an $[n+1,n+1-r]_{q}R$ code by attaching an
arbitrary column to a parity check matrix of an $[n,n-r]_{q}R$
code $\mathcal{C}$, or, equivalently, by adding an information
symbol. Clearly, by repeating the process it is possible to
obtain an $[n+\delta ,n+\delta -r]_{q}R$ code from
$\mathcal{C}$ for any integer $\delta \geq 1$. We will call
such a code a $ \delta $-extension of $\mathcal{C}$.

For a given $R\geq 1$ and for a fixed prime power $q$, let
$\mathcal{A} _{R,q} $ denote an infinite sequence of $q$-ary
linear $[n,n-r_n]_{q}R$ codes $\mathcal{C}_{n}$, $n\ge R$, with
fixed covering radius $R$. An infinite sequence
$\mathcal{A}_{R,q}$ of covering codes is called an \emph{
infinite family of covering codes }or an infinite code
family,\emph{\ }or simply infinite family.\emph{\ }

For infinite families $\mathcal{A}_{R,q}$ we consider
\emph{asymptotic covering densities}
\begin{eqnarray}
\overline{\mu }_{q}(R,\mathcal{A}_{R,q}) &=&\liminf_{n\rightarrow \infty
}\mu _{q}(n,R,\mathcal{C}_{n}).  \label{form1_liminfdens} \\
\mu _{q}^{\ast }(R,\mathcal{A}_{R,q}) &=&\limsup_{n\rightarrow \infty }\mu
_{q}(n,R,\mathcal{C}_{n}).  \label{form1_limsupdens}
\end{eqnarray}
We will write $\overline{\mu }_{q}(R)$ ($\mu _{q}^{\ast }(R)$
resp.) for $ \overline{\mu }_{q}(R,\mathcal{A}_{R,q})$ ($\mu
_{q}^{\ast }(R,\mathcal{A} _{R,q})$ resp.) if the family
$\mathcal{A}_{R,q}$ is clear from the context.

For an infinite family $\mathcal{A}_{R,q}$ the sequence of
codimensions $ r_{n}$ will be assumed to be non-decreasing. In
fact, if $r_{n+1}<r_{n}$ for some $n$, then any $1$-extension
$\mathcal{C}^{\ast }$ of $\mathcal{C}_{n}$ has a better
covering density than $\mathcal{C}_{n+1}$, and therefore it is
convenient to replace $\mathcal{C}_{n+1}$ with
$\mathcal{C}^{\ast }$.

A code $\mathcal{C}_{n}$ will be called a \emph{supporting
code} of $ \mathcal{A}_{R,q}$ if $r_{n}>r_{n-1}$, a
\emph{filling code} otherwise. It is immediately seen that a
filling code must have the same parameters of a $ \delta
$-extension of some supporting code, and this motivates our
notation. The subsequence of supporting codes will be denoted
as $\mathcal{C}_{n_{i}}$.

Throughout the paper, constructing an infinite family, we will
only describe supporting codes, whereas the filling codes will
be assumed to be obtained via $\delta$-extension. The words
\textquotedblleft to construct a family\textquotedblright\ will
mean \textquotedblleft to construct the supporting codes of a
family\textquotedblright.

In this work we will mainly deal with infinite families
$\mathcal{A}_{R,q}$ for which the lengths and the codimension
of the supporting codes $\mathcal{C }_{n_{i}}$ are linked by
some function, namely $n_{i}=f_{q}(r_{i})$ where $ f_{q}$ is an
increasing function for a fixed $q$. In most cases, an explicit
expression for the function $f_{q}$ will be given.

By (\ref{form1_covdensity}), the covering density of an
$[n+1,n+1-r]_{q}R$ is greater than that of an $[n,n-r]_{q}R$
one. Therefore,
\begin{eqnarray}
\overline{\mu }_{q}(R,\mathcal{A}_{R,q}) =\liminf_{i\rightarrow \infty
}\mu _{q}(n_{i},R,\mathcal{C}_{n_{i}}),  \label{form1_liminf_r_t} \\
\mu _{q}^{\ast }(R,\mathcal{A}_{R,q}) =\limsup_{i\rightarrow \infty }\mu
_{q}(n_{i+1}-1,R,\mathcal{C}_{n_{i+1}-1}).  \label{form1_limsup_i}
\end{eqnarray}

Note that by (\ref{form1_liminf_r_t}),(\ref{form1_limsup_i}),
the lower limit of the asymptotic covering density depends only
on the supporting codes, while the upper limit depends on
filling codes.

The size $q$ of the base field $F_{q}$ is fixed for a given
family $\mathcal{ A}_{R,q}$. But, it is natural to consider an
\emph{infinite set of families } $\mathcal{A}_{R,q}$ with fixed
$R$ and infinitely growing $q$. In most constructions,
$f_{q}(r)$ is an increasing function of $q$ for a fixed$~r.$
Therefore, a central problem for \emph{linear }covering codes
is the following:

\emph{For a fixed covering radius }$R,$ \emph{find a set of
families }$ \mathcal{A}_{R,q}$\emph{\ of }$q$\emph{-ary codes
with }$q$\emph{\ running over an infinite set of prime power,
such that the covering densities
\textrm{(\ref{form1_liminf_r_t})} and
\textrm{(\ref{form1_limsup_i})} are asymptotically as small as
possible with respect to the size of the base field~$q$}.

This problem has distinct perspectives and solutions for lower
and upper limits.

As to the lower limit \textrm{(\ref{form1_liminf_r_t})}, it can
happen that the asymptotic covering density of a family
$\mathcal{A}_{R,q}$ are bounded from above by a constant
independent of~$q$. In this case $\overline{\mu }
_{q}(R,\mathcal{A}_{R,q})=O(1)$ and the family
$\mathcal{A}_{R,q}$ is said to be \textquotedblleft
\emph{good\textquotedblright }. Accordingly, an $ [n,n-r]_{q}R$
covering code is called \textquotedblleft \emph{
short\textquotedblright } if $n=O(q^{\frac{r-R}{R}})$. By (\ref
{form1_covdensity}) and (\ref{form1_liminfdens}), a family
$\mathcal{A}_{R,q} $ consisting of short codes is good. In this
case, $f_{q}(r)=O(q^{\frac{r-R}{ R}}).$ A saturating set $K$
will be said to be \textquotedblleft \emph{small}
\textquotedblright\ if the related covering code
$\mathcal{C}_{K}$ is short.

A classical example is the direct sum \cite{Coh} of $R$ copies
of the $[ \frac{q^{i}-1}{q-1},\frac{q^{i}-1}{q-1}-i]_{q}1$
perfect Hamming codes, which gives an infinite family
$\mathcal{A}_{R,q}$ of $ [n_{i},n_{i}-r_{i}]_{q}R$ codes with
parameters
\setlength{\arraycolsep}{0.0em}
\begin{eqnarray}
\mathcal{A}_{R,q}:n_{i}=R\frac{q^{i}-1}{q-1},\text{ }r_{i}=Ri,\text{ }
i=1,2,3,\ldots ;\notag\\
\overline{\mu }_{q}(R)=O\left(\frac{R^{R}}{R!}
\right).  \label{form1_DS_Ham}
\end{eqnarray}

When the upper limit is considered, it is not possible to
obtain an upper bound independent on $q$. This depends on the
fact that
\setlength{\arraycolsep}{0.0em}
\begin{eqnarray*}
\mu_q(n_{i+1}-1,R,C_{n_{i+1}-1})=\frac{V_q(n_{i+1}-1,R)}{q^{r_i}}
=\\
\mu_q(n_{i+1},R,C_{n_{i+1}})\frac{V_q(n_{i+1}-1,R)}{V_q(n_{i+1},R)}\frac{
q^{r_{i+1}}}{q^{r_i}}.
\end{eqnarray*}
Since $r_{i+1}>r_i$, this implies that the optimal case is $\mu
_{q}^{\ast }(R,\mathcal{A}_{R,q})=O(q).$ Then the following
natural issue arises.\smallskip

\noindent \textbf{Open Problem 1. }\emph{For any covering
radius }$R\geq 2,$ \emph{construct an infinite code family
}$\mathcal{A}_{R,q}$ \emph{with} $\mu _{q}^{\ast
}(R,\mathcal{A}_{R,q})=O(q).$\smallskip

To solve Open Problem 1 it is convenient to proceed as follows.
For any given integer $\gamma $ with $0\leq \gamma \leq R-1$,
construct an infinite family $\mathcal{A}_{R,q}^{(\gamma )}$
such that its supporting codes are $ [n_{u},n_{u}-r_{u}]_{q}R$
codes with codimension $r_{u}=Ru+\gamma $ and length
$n_{u}=f_{q}^{(\gamma )}(r_{u})$, where $u\geq u_{0}$ and a
constant $ u_{0}$ may depend on the family. Considering
families of type $\mathcal{A} _{R,q}^{(\gamma )}$ is a standard
method of investigation of \emph{linear } covering codes, see
\cite{Coh},\cite{DD-L},\cite{Dav95},\cite{DavNBCR2},\cite
{DavO1,DavO3,DavO5},\cite{DMPIEEE} and the references therein;
families $ \mathcal{A}_{R,q}^{(\gamma )}$ with distinct values
of $\gamma $ often have distinct properties.

Assume that we have $R$ good infinite code families
$\mathcal{A} _{R,q}^{(\gamma )},$ $\gamma =0,1,\ldots ,R-1.$
Let us consider the infinite family
$\widehat{\mathcal{A}}_{R,q},$ whose supporting codes are the
union of those of all the families$~\mathcal{A}_{R,q}^{(\gamma
)}.$ The family $ \widehat{\mathcal{A}}_{R,q}$ contains an
infinite sequence of $ [n_{j},n_{j}-j]_{q}R$ codes
$\mathcal{C}_{j}$ with length $ n_{j}=f_{q}^{(\gamma
_{j})}(j)$, $\gamma _{j}\equiv j\pmod R$, where $j\geq j_{0}$
and $j_{0}$ is a constant depending of constants $u_{0}$ of the
starting families. Note that it may occur that $n_{v+1} \leq
n_{v}$ for some $v$. In this case we replace the code
$\mathcal{C}_{v}$ by an $[n_{v+1}-1,n_{v+1}-1-v]_{q}R$ code
that always can be obtained from $\mathcal{C}_{v+1}$ by
removing a redundancy symbol and a suitable parity check.
Arguing as before,
\begin{equation*}
\mu _{q}^{\ast }(R,\widehat{\mathcal{A}}_{R,q})=\limsup_{j\rightarrow \infty
}\frac{V_{q}(n_{j+1},R)}{q^{j+1}}\frac{V_{q}(n_{j+1}-1,R)}{V_{q}(n_{j+1},R)}
\frac{q^{j+1}}{q^{j}}.
\end{equation*}
Since all families $\mathcal{A}_{R,q}^{(\gamma )}$ are good, we
have $ V_{q}(n_{j+1},R)/q^{j+1}=O(1)$. Hence,
\begin{equation*}
\mu _{q}^{\ast }(R,\widehat{\mathcal{A}}_{R,q})=O(q).
\end{equation*}

So, to solve Open Problem 1 it is sufficient to find a solution
to Open Problem 2.\smallskip

\noindent \textbf{Open Problem 2. } \emph{For any covering
radius $R\geq 2,$ construct $R$ infinite code families
$\mathcal{A}_{R,q}^{(0)},\mathcal{A}_{R,q}^{(1)},\ldots
,\mathcal{A}_{R,q}^{(R-1)}$ such that for each $\gamma
=0,1,\ldots ,R-1$ the supporting codes of
$\mathcal{A}_{R,q}^{(\gamma )}$ are $[n_{u},n_{u}-r_{u}]_{q}R$
codes with codimension $r_{u}=Ru+\gamma $ and length
$n_{u}=f_{q}^{(\gamma )}(r_{u})$ with $f_{q}^{(\gamma
)}(r)=O(q^{\frac{r-R}{R}})$ and }$u\geq u_{0} $\emph{\ where a
constant }$u_{0}$\emph{\ may depend on the family.}\smallskip

On one hand, infinite families $\mathcal{A}_{R,q}^{(0)}$ are
provided by example~(\ref {form1_DS_Ham}); for $R=2,3,$
families $\mathcal{A}_{R,q}^{(0)}$ with better parameters are
obtained in \cite{Dav95},\cite{DavNBCR2},\cite{DavO5}. On the
other hand, for $\gamma \geq 1$, code families
$\mathcal{A}_{R}^{(\gamma )}$ with density $\overline{\mu
}_{q}(R,\mathcal{A}_{R}^{(\gamma )})=O(1)$ are only known for
$R=2,$ $\gamma =1,$ $q=(q^{\prime })^{2}$ \cite{DavNBCR2}, and
$R=3,$ $\gamma =1,$ $q=(q^{\prime })^{3}$~\cite{DMPIEEE}.

In this paper, Open Problem 2 (and Open Problem 1) is solved
for an arbitrary covering radius $R\geq 2$ and $q=(q^{\prime
})^{R}$ where $ q^{\prime }$ is a power of prime.

Our main tools are the $q^{m}$-concatenating constructions of covering
codes, and the connection between covering codes and saturating sets in
projective spaces.

The $q^{m}$-concatenating constructions are proposed in
\cite{DavPPI} and are developed in
\cite{DavParis},\cite[Supplement]
{Rene-PhD},\cite{DD-L,Dav95,davShumen,Dav-Nonlin},\cite{DavNBCR2},\cite{DavO1,DavO3,DavO5,Dav-2001-newconst},\cite
{DFMP-IEEE-LO}, see also \cite[Sec.\thinspace 5.4]{Coh} and
\cite {Handbook-coverings}. These constructions are the
fundamental instrument for obtaining infinite families of
covering codes with a fixed radius. Using a starting code as a
\textquotedblleft seed\textquotedblright , the $q^{m}$
-concatenating constructions yield an infinite family of new
codes with the same covering radius and with almost the same
covering density. If the starting code is short then the new
infinite family is good.

Linear codes arising from small saturating sets are a
convenient choice for the starting codes of the
$q^{m}$-concatenating constructions \cite{Dav95},
\cite{DavNBCR2},\cite{DavO5},\cite{DMPIEEE},\cite{DFMP-IEEE-LO}.

The achievements of the present paper are mainly a consequence
of new constructions of small\ saturating sets, some of which
rely on the concept of a \emph{multifold strong blocking set}
that is introduced in this work. We have also thoroughly
analyzed and collected the known results on the upper bounds on
the length function, in particular for the cases $R=2,3$. We
have updated tables about the upper bounds and formulas for
infinite code families. As a result of our previously mentioned
constructions, many new upper bounds on the length function are
obtained.

The paper is organizes as follows. In Section
\ref{Sec2_q^m_concat} the $ q^{m}$-concatenating constructions,
used in this work, are recalled. In Section \ref{Sec4_NewSmSat}
new constructions of small $\varrho $ -saturating sets,
including those relying on the new concept of strong blocking
sets, are described. Section \ref{Sec5_Tab} contains updated
tables about the upper bounds on $\ell _{q}(r,R)$ for $R=2,3,$
$r=3,4,5.$ In Sections \ref{Sec6_R=2},\ref{Sec7_R=3}, and
\ref{Sec8_R>=4} we consider codes with covering radii $R=2,$
$R=3,$ and $R\geq 4$. Section \ref{Sec9_R nonprime} provides
results for nonprime covering radius.

Some of the results from this work were briefly presented without proofs in
\cite{DGMP-ACCT2008},\cite{DGMP-Petersb2008}.

\section{\label{Sec2_q^m_concat}$q^{m}$-Concatenating Constructions
Lengthening\newline
Covering Codes}

In this section we describe the common ideas and the popular
versions of the $q^{m}$-concatenating constructions. Other
versions can be found in \cite[ Sec.\thinspace
5.4]{Coh},\cite{Handbook-coverings},\cite{DavPPI},\cite
{DavParis},\cite[Supplement]{Rene-PhD},\cite{DD-L,Dav95,davShumen,Dav-Nonlin},\cite
{DavO5},\cite{Dav-2001-newconst},\cite{DMPIEEE},\cite{DFMP-IEEE-LO}.
Specific constructions for $R=2$ are given in detail
in~\cite{DavNBCR2}.

Using a starting $[n_{0},n_{0}-r_{0}]_{q}R$ code of length
$n_{0},$ the $ q^{m}$-concatenating constructions yield an
infinite family of $ [n,n-(r_{0}+Rm)]_{q}R$ codes with the same
covering radius $R$ and length $ n=q^{m}n_{0}+N_{m}$, where $m$
ranges over an infinite set of integers. Here $N_{m}\leq
R\theta _{m,q}$, where
\begin{equation*}
\theta _{m,q}=\frac{q^{m}-1}{q-1}.
\end{equation*}
It should be noted that all $q^{m}$-concatenating constructions
have the contribution $q^{m}n_{0}$ into $n$; two of them may
differ by the value of $ N_{m}.$

Throughout this paper, all matrices and columns are $q$-ary. An
element of $ F_{q^{m}}$ written in a $q$-ary matrix denotes an
$m$-dimensional column containing its coordinates with respect
to a fixed basis of $F_{q^{m}}$ over $F_{q}$; viceversa, an
$m$-dimensional vector can be viewed as an element of$
~F_{q^{m}}$.

\subsection{$(R,\ell) $-partitions and $(R,\ell) $
-objects}

\begin{definition}
\label{Def2_R,l code}Let $\mathbf{H}$ be a parity-check matrix
of an $ [n,n-r]_{q}R$ code $V$ and let $0\leq \ell \leq R$.

\emph{i) \ }A partition of the column set of the matrix
${\mathbf{H}}$ into nonempty subsets is called an \linebreak
$(R,\ell )$\emph{-partition} if every column of $F_{q}^{r}$
(including the zero column) is equal to a linear combination
with nonzero coefficients of at least $\ell $ and at most $R$
columns of $\mathbf{H}$ belonging to distinct subsets. For an
$(R,0)$ -partition we can formally treat the zero column as the
linear combination of 0 columns.

An $R$-partition is an $(R,\ell)$-partition for some $\ell\ge 0$.

\emph{ii) }If $\mathbf{H}$ admits an $(R,\ell) $-partition, the
code $V$ is called an $(R,\ell) $\emph{-object }and is denoted
as \linebreak an $ [n,n-r]_{q}R,\ell $ code or an
$[n,n-r,d]_{q}R,\ell $ code, where $d$ is the minimum distance
of $V$.
\end{definition}

Clearly, the \emph{trivial} partition of a parity-check matrix
of an $ [n,n-r]_{q}R,\ell $ code into $n$ one-element subsets
is an $(R,\ell )$ -partition.

Note that in Definition \ref{Def2_R,l code}, it is not
necessary that $\ell $ is the greatest value with the
properties considered. Any $(R,\ell) $ -partition with $\ell
>0$ is also an $(R,\ell _{1})$-partition with $\ell
_{1}=0,1,\ldots ,\ell -1$.

\begin{lemma}
\label{Lem2_d-and-l}\cite{DavPPI},\cite{DD-L},\cite{Dav95} An $
[n,n-r,d]_{q}R $ code is an $[n,n-r,d]_{q}R,\ell $ code with
$\ell \geq 1$ if and only if $d\leq R.$ If $d>R$ the maximum
possible value of $\ell $ is zero$.$
\end{lemma}

A \emph{spherical }$(R,\ell) $-\emph{capsule }with center $c$
in $F_{q}^{n}$ is the set $\{v:v\in F_{n}^{q},$ $0\leq \ell
\leq d(v,c)\leq R\}$ (see \cite {DavPPI}). It is easy to see
that spherical $(R,\ell) $-capsules centered at vectors of an
$(R,\ell) $-object cover the space $F_{q}^{n}.$

\subsection{ Basic $q^{m}$-Concatenating Constructions}

We give a basic $q^{m}$-concatenating construction QM\textbf{\ }based on
ideas in~\cite{DavPPI},\cite{DD-L},\cite{Dav95},\cite{Dav-Nonlin}.\medskip

\noindent \textbf{Basic Construction QM.} Let
$\mathbf{H}_{0}=[\mathbf{h}_{1} \mathbf{h}_{2}\mathbf{\ldots
h}_{n_{0}}]$, with $\mathbf{h}_{j}\in F_{q}^{r_{0}}$, be a
parity check matrix of an $[n_{0},n_{0}-r_{0}]_{q}R, \ell _{0}$
\emph{starting} code $V_{0}$. Assume that $\mathbf{H}_{0}$ has
a starting $(R,\ell _{0})$-partition $\mathcal{P}_{0}$ into
$p_{0}$ subsets. Let $m\geq 1$ be an integer parameter
depending on $p_{0}$ and $n_{0}$. To each column
$\mathbf{h}_{j}$ we associate an element $\beta _{j}$ $\in
F_{q^{m}}\cup \{\ast \}$ so that $\beta _{i}\neq \beta _{j}$ if
columns $ \mathbf{h}_{i}$ and $\mathbf{h}_{j}$ belong to
\emph{distinct }subsets of~$ \mathcal{P}_{0}.$ If
$\mathbf{h}_{i}$ and $\mathbf{h}_{j}$ belong to the same subset
we are free to assign either $\beta _{i}=\beta _{j}$ or $\beta
_{i}\neq \beta _{j}$. We call $\beta _{j}$ an \emph{indicator}
of column $ \mathbf{h}_{j}.$ Let $\mathcal{B}=\{\beta
_{1},\beta _{2},\ldots ,\beta _{n_{0}}\}$ be an \emph{indicator
set}. It is necessary that $|\mathcal{B} |\geq p_{0}.$ Also,
let $\mathbf{C}$ be an $(r_{0}+Rm)\times N_{m}$ matrix with
$N_{m}\leq (R-\ell _{0})\theta _{m,q}$. Finally, define $V$ as
the $ [n,n-(r_{0}+Rm)]_{q}R_{V}$ code with $n=q^{m}n_{0}+N_{m}$
and the parity-check matrix of the form
\begin{eqnarray}
\mathbf{H}_{V}&\,=\,&\left[ \mathbf{C}~\mathbf{B}_{1}~\mathbf{B}_{2}~\ldots ~
\mathbf{B}_{n_{0}}\right] ,  \label{form2_QM-H}\\
\mathbf{B}_{j}&\,=\,&\left[ \renewcommand{\arraystretch}{1.1}
\begin{array}{@{}c@{\,\,\,}c@{\,}c@{\,}c}
\mathbf{h}_{j} & \mathbf{h}_{j} & \mathbf{\cdots } & \mathbf{h}_{j} \\
\xi _{1} & \xi _{2} & \cdots & \xi _{q^{m}} \\
\beta _{j}\xi _{1} & \beta _{j}\xi _{2} & \cdots & \beta _{j}\xi _{q^{m}} \\
\beta _{j}^{2}\xi _{1} & \beta _{j}^{2}\xi _{2} & \cdots & \beta _{j}^{2}\xi
_{q^{m}} \\
\vdots & \vdots & \vdots & \vdots \\
\beta _{j}^{R-1}\xi _{1} & \beta _{j}^{R-1}\xi _{2} & \cdots & \beta
_{j}^{R-1}\xi _{q^{m}}
\end{array}
\right]\text{if }\beta _{j}\in F_{q^{m}},\notag\\
\mathbf{B}_{j}&\,=\,&\left[
\renewcommand{\arraystretch}{1.0}
\begin{array}{@{\,}c@{\,\,}c@{\,}c@{\,}c}
\mathbf{h}_{j} & \mathbf{h}_{j} & \mathbf{\cdots } & \mathbf{h}_{j} \\
0 & 0 & \cdots & 0 \\
\vdots & \vdots & \vdots & \vdots \\
0 & 0 & \cdots & 0 \\
\xi _{1} & \xi _{2} & \cdots & \xi _{q^{m}}
\end{array}
\right] \text{ if }\beta _{j}=\ast ,\notag
\end{eqnarray}
where $\{\xi _{1},\xi _{2},\ldots ,\xi _{q^{m}}\}=F_{q^{m}},$
$\xi _{1}=0,$ $ \xi _{2}=1$. Note that the submatrix
$\mathbf{C}$ is not needed if $\ell _{0}=R$.

If $m,\mathbf{C}$ and $\mathcal{B}$ are carefully chosen, then
the covering radius $R_{V}$ of the new code $V$ is equal to the
covering radius $R$ of the starting code $V_{0}.$ Examples are
shown in Constructions QM$_{1}$ - QM$ _{8}$ below.

We use the following notations: \newline $\mathbf{W}_{m}$ is a
parity-check matrix of the $[\theta _{m,q},\theta
_{m,q}-m]_{q}1$ Hamming code$;$ \newline $\mathbf{A}_{R^{\prime
},m}$ is a parity-check matrix of an $[n^{\prime },n^{\prime
}-R^{\prime }m]_{q}R^{\prime }$ code $V_{R^{\prime },m}$ (in
most cases we will assume that either $n^{\prime
}=\overline{\ell } _{q}(R^{\prime }m,R^{\prime })$ or
$n^{\prime }=\ell _{q}(R^{\prime }m,R^{\prime })$); \newline
$\mathbf{0}_{k}$ is the zero matrix with $k$ rows (the number
of columns will be clear by context); \newline $\mathbf{\Sigma
}_{R^{\prime \prime },m}$ is the \textquotedblleft direct
sum\textquotedblright\ of $R^{\prime \prime }$ matrices
$\mathbf{W}_{m},$ i.e. an $R^{\prime \prime }m\times R^{\prime
\prime }\theta _{m,q}$ matrix of the form
\begin{equation}
\mathbf{\Sigma }_{R^{\prime \prime },m}\mathbf{=}\left[
\renewcommand{\arraystretch}{1.0}
\begin{array}{@{\,\,}c@{\,\,\,}c@{\,\,}c@{\,\,}c@{\,\,}}
\mathbf{W}_{m} & \mathbf{0}_{m} & \cdots & \mathbf{0}_{m} \\
\mathbf{0}_{m} & \mathbf{W}_{m} & \cdots & \mathbf{0}_{m} \\
\vdots & \vdots & \ddots & \vdots \\
\mathbf{0}_{m} & \mathbf{0}_{m} & \cdots & \mathbf{W}_{m}
\end{array}
\right] .  \label{form2_Sum_R'',m}
\end{equation}
Note that $\mathbf{\Sigma }_{R^{\prime \prime },m}$ is a
parity-check matrix of an $[R^{\prime \prime }\theta
_{m,q},R^{\prime \prime }\theta _{m,q}-R^{\prime \prime
}m]_{q}R^{\prime \prime }$ code, see the direct sum
construction in Section \ref{Sec6_R=2}.

\noindent \textbf{Construction QM}$_{1}$\textbf{.} Here $R\geq 2,$ $\ell
_{0}=0,$ $q^{m}+1\geq p_{0}$, $\mathcal{B}\subseteq F_{q^{m}}\cup \{\ast \},$
\begin{equation}
\mathbf{C=}\left[ \renewcommand{\arraystretch}{1.0}
\begin{array}{@{\,\,}c@{\,\,}}
\mathbf{0}_{r_{0}} \\
\mathbf{\Sigma }_{R,m}
\end{array}
\right] ,\text{ }n=q^{m}n_{0}+R\theta _{m,q}.  \label{form2_QM_1}
\end{equation}
\noindent \textbf{Construction QM}$_{2}$\textbf{.} Here $R\geq 2,$ $1\leq
\ell _{0}<R,$ $q^{m}\geq p_{0}$, $\mathcal{B}\subseteq F_{q^{m}},$
\begin{equation}
\mathbf{C=}\left[ \renewcommand{\arraystretch}{1.0}
\begin{array}{@{\,\,}c@{\,\,}}
\mathbf{0}_{r_{0}+\ell _{0}m} \\
\mathbf{\Sigma }_{R-\ell _{0},m}
\end{array}
\right] ,\text{ }n=q^{m}n_{0}+(R-\ell _{0})\theta _{m,q}.  \label{form2_QM_2}
\end{equation}

\noindent \textbf{Construction QM}$_{3}$\textbf{.} Here $R\geq 2,$ $\ell
_{0}=R,$ $q^{m}+1\geq p_{0}$, $\mathcal{B}\subseteq F_{q^{m}}\cup \{\ast \},$
\begin{equation}
\mathbf{C}\text{ is absent},\text{ }n=q^{m}n_{0}.  \label{form2_QM_3}
\end{equation}

\noindent \textbf{Construction QM}$_{4}$\textbf{.} Here $R\geq
2,$ $\ell _{0}=0,$ $q^{m}-1\geq p_{0}$, $\mathcal{B}\subseteq
F_{q^{m}}^{\ast },$ \setlength{\arraycolsep}{0.0em}
\begin{eqnarray}
\mathbf{C}&\,=\,&\left[ \renewcommand{\arraystretch}{1.0}
\begin{array}{@{\,\,}c@{\,\,}c@{\,\,}}
\mathbf{0}_{r_{0}} & \mathbf{0}_{r_{0}} \\
\mathbf{\Sigma }_{\left\lfloor R/2\right\rfloor ,m} & \mathbf{0}
_{\left\lfloor R/2\right\rfloor m} \\
\mathbf{0}_{\left\lceil R/2\right\rceil m} & \mathbf{A}_{\left\lceil
R/2\right\rceil ,m}
\end{array}
\right] ,\notag\\
n&\,\leq\,& q^{m}n_{0}+\left\lfloor R/2\right\rfloor
\theta _{m,q}+\overline{\ell }_{q}(\left\lceil R/2\right\rceil
m,\left\lceil R/2\right\rceil ).  \label{form2_QM_4}
\end{eqnarray}

\subsection{$q^{m}$-Concatenating Constructions with a
Complete Set of Indicators (CSI)}

In these versions of the basic Construction QM we \emph{must
}use \emph{all} elements of $F_{q^{m}}$ or $F_{q^{m}}\cup
\{\ast \}$ as indicators $\beta _{j}.$ To this end, perhaps, we
should assign distinct indicators to columns from the same
subset of an $R$-partition. As a result the size of the
submatrix $\mathbf{C}$ is reduced. \medskip

\noindent \textbf{Construction QM}$_{5}$\textbf{.} Here $R\geq 2,$ $\ell
_{0}=0,$ $n_{0}\geq q^{m}\geq p_{0},$ $\mathcal{B}=F_{q^{m}},$
\begin{equation}
\mathbf{C=}\left[ \renewcommand{\arraystretch}{0.9}
\begin{array}{@{\,\,}c@{\,\,}}
\mathbf{0}_{r_{0}+m} \\
\mathbf{\Sigma }_{R-1,m}
\end{array}
\right] ,\text{ }n=q^{m}n_{0}+(R-1)\theta _{m,q}.  \label{form2_QM_5}
\end{equation}

\noindent \textbf{Construction QM}$_{6}$\textbf{.} Here $R\geq
3,$ $\ell _{0}=0,$ $n_{0}\geq q^{m}+1\geq p_{0},$
$\mathcal{B}=F_{q^{m}}\cup \{\ast \}, $
\begin{equation}
\mathbf{C=}\left[ \renewcommand{\arraystretch}{0.9}
\begin{array}{@{\,\,}c@{\,\,}}
\mathbf{0}_{r_{0}+m} \\
\mathbf{\Sigma }_{R-1,m}
\end{array}
\right] ,\text{ }n=q^{m}n_{0}+(R-1)\theta _{m,q}.  \label{form2_QM6}
\end{equation}
\noindent \textbf{Construction QM}$_{7}$\textbf{.} Here $R=3,$
$\ell _{0}=0,$ $n_{0}\geq q^{m}+1\geq p_{0}$,
$\mathcal{B}=F_{q^{m}}\cup \{\ast \},$ $ q=2^{i},$
\begin{equation}
\mathbf{C=}\left[ \renewcommand{\arraystretch}{0.9}
\begin{array}{@{\,\,}c@{\,\,}}
\mathbf{0}_{r_{0}+m} \\
\mathbf{W}_{m} \\
\mathbf{0}_{m}
\end{array}
\right] ,\text{ }n=q^{m}n_{0}+\theta _{m,q}.  \label{form2_QM_6}
\end{equation}
\noindent \textbf{Construction QM}$_{8}$\textbf{.} Here $R=4,$
$\ell _{0}=0,$ $n_{0}\geq q^{m}\geq p_{0}$,
$\mathcal{B}=F_{q^{m}},$ $3$ does not divide $ q^{m}-1,$
\begin{equation}
\mathbf{C=}\left[ \renewcommand{\arraystretch}{0.9}
\begin{array}{@{\,\,}c@{\,\,}c@{\,\,}}
\mathbf{0}_{r_{0}+m} & \mathbf{0}_{r_{0}+m} \\
\mathbf{A}_{2,m} & \mathbf{0}_{2m} \\
\mathbf{0}_{m} & \mathbf{W}_{m}
\end{array}
\right] ,\text{ }n\leq q^{m}n_{0}+\overline{\ell }_{q}(2m,2)+\theta _{m,q}.
\label{form2_QM_7}
\end{equation}

Other constructions CSI including these with $n_{0}<q^{m}$ can
be found in \cite{Dav95},\cite{Dav-2001-newconst}.

\subsection{ Summary}

\begin{theorem}
\label{Th2_basQM}\cite{Coh},\cite{DavPPI},\cite{Dav95},\cite{DavNBCR2},\cite
{DMPIEEE},\cite{DavO5} In all Constructions QM$_{i}$ the new
code $V$ is an $ [n,n-(r_{0}+Rm),3]_{q}R,\ell $ code with
covering radius $R$ and $\ell \geq \ell _{0}$.
\end{theorem}

\begin{corollary}
It holds that
\begin{equation}
\ell _{q}(r_{0}+Rm,R)\leq q^{m}\ell _{q}(r_{0},R)+\label{form2_QM1-2}\\
\end{equation}
\begin{eqnarray*}
\left\{
\begin{array}{ll}
R\theta _{m,q}\medskip & \text{ if }q^{m}+1\geq \ell _{q}(r_{0},R) \\
\left\lfloor \frac{R}{2}\right\rfloor \theta _{m,q}+\ell _{q}(\left\lceil
\frac{R}{2}\right\rceil m,\left\lceil \frac{R}{2}\right\rceil ) & \text{ if }
q^{m}>\ell _{q}(r_{0},R)
\end{array}
\right. .
\end{eqnarray*}
\end{corollary}

\begin{IEEEproof}
In Constructions QM$_{1}$ and QM$_{4}$ we put $n_{0}=\ell
_{q}(r_{0},R)$ and then use the trivial $R$-partition. In the
code $V_{R^{\prime },m}$ with $ R^{\prime }=\left\lceil
\frac{R}{2}\right\rceil $ we put $n^{\prime }=\ell
_{q}(\left\lceil \frac{R}{2}\right\rceil m,\left\lceil
\frac{R}{2} \right\rceil )$.
\end{IEEEproof}

Note that Constructions QM$_{1}$-QM$_{4}$ provide an
\emph{infinite family } of the $\emph{new}$
$[n,n-(r_{0}+Rm)]_{q}R$ codes $V$ with growing codimension
$r=r_{0}+Rm.$ In Constructions QM$_{5}$-QM$_{8}$ instead, the
value of $m$ cannot assume arbitrarily large values. However,
these construction can be used in an iterative process where
the new codes are the starting ones for the following steps
\cite{Dav95},\cite{DavNBCR2}. As result we obtain an
\emph{infinite code family, }see, e.g., \cite[ Rem.\thinspace
5]{DavNBCR2}.\emph{\ }For this iterative process, it is
important that in the new codes obtained by the
$q^{m}$-concatenating constructions the value of $\ell $ is
increasing and eventually reaches $R,$ see \cite[Sec.\thinspace
IV]{Dav95} and Examples in Section~\ref{Sec7_R=3}.

\begin{remark}
\emph{i) }By
(\ref{form1_covdensity}),(\ref{form1_liminfdens}),(\ref
{form1_liminf_r_t}),(\ref{form2_QM_1}), in Construction
QM$_{1}$ the covering density of the starting code $V_{0}$ and
the lower limit of the asymptotic covering density of the
infinite family of the new codes $V$ are, respectively, $\mu
_{q}(n_{0},R)\approx \frac{(q-1)^{R}n_{0}^{R}}{R!q^{r_{0}} }$
and $\overline{\mu }_{q}(R)\approx
\frac{(q-1)^{R}(n_{0}+R/q)^{R}}{ R!q^{r_{0}}}.$ We have
$\overline{\mu }_{q}(R)/\mu _{q}(n_{0},R)\approx (1+
\frac{R}{qn_{0}})^{R}.$ This shows that for the
$q^{m}$-concatenating constructions the lower limit of the
asymptotic covering density of the new family is somewhat
greater than covering density of the starting code. However, it
should be noted that the difference is not significant when the
value of $R/qn_{0}$ is small.

\emph{ii) }By (\ref{form2_QM_1})-(\ref{form2_QM_7}), if the
starting code $ V_{0}$ is \textquotedblleft
short\textquotedblright , i.e. $n_{0}=O(q^{\frac{
r_{0}-R}{R}}),$ then the all new $[n,n-(r_{0}+Rm)]_{q}R$ codes
$V,$ obtained by the $q^{m}$-concatenating constructions, are
\textquotedblleft short\textquotedblright\ too, i.e.
$n=O(q^{\frac{r_{0}+Rm-R}{R}}).$ This means that the infinite
family of the codes $V$ is \textquotedblleft
good\textquotedblright\ with $\overline{\mu }_{q}(R)=O(1).$
\end{remark}

\section{\label{Sec4_NewSmSat}New \textquotedblleft Small\textquotedblright\
Saturating Sets}

\subsection{Multifold Strong Blocking Sets}

In a projective space a $t$-fold blocking set with respect to
subspaces of some fixed dimension is a set of points that meets
every such subspace in at least $t$ points. To describe new
constructions of relatively small $\rho $ -saturating sets in
spaces $PG(v,q)$ with $q=(q^{\prime })^{\rho +1}$ we introduce
a new concept of $t$-fold \textbf{strong }blocking set.

\begin{definition}
\label{Def4_strong_block} Let $2\leq t\leq v.$ A pointset $B$
in a projective space $PG(v,q)$ is a $\emph{t}$\emph{-fold
strong blocking set} if every $(t-1)$-dimensional subspace of
$PG(v,q)$ is spanned by $t$ points in $B$.
\end{definition}

Let $(x_{0},x_{1},\ldots ,x_{v}),$ where $x_{i}\in F_{q},$ be
homogenous coordinates for a point $P$ in $PG(v,q)$ and let $
P^{u}=(x_{0}^{u},x_{1}^{u},\ldots ,x_{v}^{u}).$

\begin{theorem}
\label{Th4_quattro} Let $\rho $ be any positive integer. Let $q=(q^{\prime
})^{\rho +1}$. Let $v\geq \rho +1$. Any $(\rho +1)$-fold strong blocking set
in a subgeometry $PG(v,q^{\prime })\subset PG(v,q)$ is a $\rho $-saturating
set in the space $PG(v,q)$.
\end{theorem}

\begin{IEEEproof}
Let $B$ be a $(\rho +1)$-fold strong blocking set in
$PG(v,q^{\prime })$. Let $P$ be a point in $PG(v,q)\setminus
B$. By definition of $(\rho +1)$ -fold strong blocking set we
only need to show that there exists a $\rho $ -dimensional
subspace of $PG(v,q^{\prime })$ passing through $P$. Consider
the subspace $\Sigma (P)$ of $PG(v,q)$ generated by the point
set $ O(P):=\{P,P^{q^{\prime }},P^{(q^{\prime })^{2}},\ldots
,P^{(q^{\prime })^{\rho }}\}$. As $\left( {P^{(q^{\prime
})^{\rho }}}\right) ^{q^{\prime }}=P^{q}=P,$ the Frobenius
collineation $X\mapsto X^{q^{\prime }}$ fixes $ O(P)$.
Therefore $\Sigma (P)$ is a subspace of $PG(v,q^{\prime })$.
Clearly $ \Sigma (P)$ is contained in some $\rho $-dimensional
subspace of $ PG(v,q^{\prime })$ (if the points in $O(P)$ are
independent, then this subspace coincides with $\Sigma (P) $).
As $P\in \Sigma (P)$, the assertion is proved.
\end{IEEEproof}

\subsection{Small $\rho $-Saturating Sets in Spaces
 $PG(\rho +1,(q^{\prime })^{\rho +1})$}

\begin{corollary}
\label{Cor4_2-fold} Let $q=(q^{\prime })^{2}.$ Any $2$-fold blocking set in
the subplane $PG(2,\sqrt{q})\subset PG(2,q)$ is a $1$-saturating set in the
plane $PG(2,q)$.
\end{corollary}

\begin{IEEEproof}
As a line is spanned by any two its points, a 2-fold blocking
set in a projective plane is always a 2-fold \emph{strong
}blocking set. Then we use Theorem~\ref{Th4_quattro}.
\end{IEEEproof}

Note that Corollary \ref{Cor4_2-fold} is also given in \cite{KKKPS}.

\begin{theorem}
\label{Th4_q4-1sat} Let $q=(q^{\prime })^{4}.$ In $PG(2,q)$
there is a $1$ -saturating set of size
$2\sqrt{q}+2\sqrt[4]{q}+2$.
\end{theorem}

\begin{IEEEproof}
The union of two disjoint Baer subplanes in $PG(2,\sqrt{q})$ is
a 2-fold blocking set \cite{balS96}. Then we use Corollary
\ref{Cor4_2-fold}.
\end{IEEEproof}

In $PG(2,q),$ $q$ not a square, 2-fold blocking sets of size $b\leq 3q-2$
are not known in the literature \cite{balS96},\cite{BallHirs}. We give here
some results for $q=p^{3},$ $p$ prime.

\begin{theorem}
\label{Th4_q3_2fold}Let $q=p^{3},$ $p$ prime, $p\leq 73$. Then
in $PG(2,q)$ there is a $2$-fold blocking set of size $2\left(
q+\sqrt[3]{q^{2}}+\sqrt[3]{ q}+1\right) .$
\end{theorem}

\begin{IEEEproof}
By \cite[Lem. 13.8 (iii)]{Hirs}, the point set
\begin{eqnarray*}
B=\{(1,x,x^{p})\mid x\in F_{q}\}\cup \\
\{(0,1,m)\mid m\in F_{q},
m^{p^{2}+p+1}=1\}
\end{eqnarray*}
is a $1$-blocking set in $PG(2,q)$ of size $q+p^{2}+p+1.$ We are looking for
a projectivity $\gamma $ for which $B\cap \gamma (B)=\emptyset $ holds. Then
$B\cup \gamma (B)$ is a $2$-fold blocking set in $PG(2,q).$

Let $H$ be the multiplicative subgroup of $F_{q}^{\ast }$
consisting of the $ (p-1)$th powers in $F_{q}$ (equivalently,
$H=\{y\in F_{q}\mid y^{p^{2}+p+1}=1\}$). For $a,b\in
F_{q}^{\ast }$, $b\notin H,$ we consider the projectivity
$\gamma _{a,b}(r,s,t)=(t-r,abr,as) $. Obviously, $\gamma
_{a,b}(0,1,m)=(m,0,a)=(1,0,a/m)\notin B.$ Also, $\gamma
_{a,b}(1,1,1)=(0,ab,a)=(0,b,1)\notin B$ as $b^{p^{2}+p+1}\neq
1$. Finally, for $x\neq 1,$ $\gamma
_{a,b}(1,x,x^{p})=(1,ab/(x^{p}-1),ax/(x^{p}-1))\in B$ if and
only if $a^{p-1}b^{p}=(x^{p}-1)^{p-1}x.$

So, $B\cap \gamma _{a,b}(B)=\emptyset $ if and only if the
equation $ a^{p-1}b^{p}=(x^{p}-1)^{p-1}x$ has no solution
in$~F_{q}$.

Now, note that any element $c\in F_{q}^{\ast }\setminus H$ can be expressed
as a product $a^{p-1}b^{p}$ with $a,b\in F_{q}^{\ast }$, $b\notin H$. In
fact, $c$ belongs to some coset $dH$, $d\notin H$, and therefore $c=da^{p-1}$
for some $a\in F_{q}^{\ast }$. Let $b=d^{p^{2}}\notin H$, so that $b^{p}=d$.

Then the following claim is proved: if there is an element $c\in F_{q}^{\ast
}\setminus H$ such that $c\neq (x^{p}-1)^{p-1}x$ for any $x\in F_{q}$, then
there exist $a,b$ such that $B\cap \gamma _{a,b}(B)=\emptyset $.

The existence of such element $c$ has been tested by computer
for every prime $p\leq 73.$
\end{IEEEproof}

\begin{corollary}
\label{Cor4_q6-1sat} Let $q=(q^{\prime })^{6},$ $q^{\prime }$
prime, $ q^{\prime }\leq 73$. In $PG(2,q)$ there is a
$1$-saturating set of size $2
\sqrt{q}+2\sqrt[3]{q}+2\sqrt[6]{q}+2$.
\end{corollary}

Note that the smallest previously known 1-saturating sets in
$PG(2,q),$ $ q=(q^{\prime })^{2},$ have size $3\sqrt{q}-1$
\cite[Th.\ 5.2]{Dav95}, cf. Theorem \ref{Th4_q4-1sat} and
Corollary \ref{Cor4_q6-1sat}.

Now we construct a 3-fold strong blocking set in $PG(3,q)$. Let
$ l_{1},l_{2},l_{3}$ be the lines with the following equations:
\begin{equation*}
l_{1}:x_{0}=x_{2}=0;\text{ }l_{2}:x_{1}=x_{3}=0;\text{ }l_{3}:x_{0}=x_{3},
\text{ }x_{1}=x_{2}.
\end{equation*}
These lines are pairwise skew, and are all contained in the hyperbolic
quadric $\mathcal{Q}:x_{0}x_{1}=x_{2}x_{3}$. Let $g$ be any line disjoint
from $\mathcal{Q}$, and let
\begin{equation}
B=l_{1}\cup l_{2}\cup l_{3}\cup g.  \label{form4_R=3_r=3t+1_2sat-4lines}
\end{equation}
A possible choice for $g$ is the following:
\begin{eqnarray*}
g:\left\{
\begin{array}{l}
x_{0}=x_{1},\,\,  x_{2}=kx_{3},\,\,   k\text{ non-square in }F_{q}, \,\,
 \text{if }q \text{ odd}. \\
x_{0}=x_{1}+x_{3},\,\,   x_{2}=kx_{3},\,\, \\\phantom{x_{0}=x_{1}+} T^{2}+T+k\text{ irreducible
over }F_{q},\,\, \text{if }q\text{ even.}
\end{array}
\right. .
\end{eqnarray*}

\begin{theorem}
\label{Th4_dim6} The set $B$ of $(\ref{form4_R=3_r=3t+1_2sat-4lines})$ has
size $4q+4$ and it is a $3$-fold strong blocking set in $PG(3,q)$.
\end{theorem}

\begin{IEEEproof}
We need to show that any plane $\pi $ of $PG(3,q)$ meets $B$ in
three non collinear points. If one of lines of $B$ lies on $\pi
$, then the assertion is trivial. Let $P_{i}=\pi \cap l_{i}$.
Assume that $P_{1},P_{2},P_{3}$ are collinear. Then the line
$l$ through $P_{1},P_{2}$ and $P_{3}$ is contained in
$\mathcal{Q}$, by the\textquotedblleft \emph{three then
all\textquotedblright } principle for quadrics in projective
spaces. As $ R=\pi \cap g\notin \mathcal{Q}$, we have that $R$
is not collinear with $ P_{1}$ and $P_{2}$.
\end{IEEEproof}

\begin{remark}
Any $3$-fold strong blocking set $B$ in $PG(3,q)$ has at least $3q+3$
points. Let $l$ be any line such that $l\cap B=\emptyset $. Then each of the
$q+1$ planes in the pencil through $l$ must contain three points of $B$.
\end{remark}

\begin{corollary}
\label{Cor4_3fold_PG(3,q)} Let $q=(q^{\prime })^{3}$. In $PG(3,q)$ there is
a $2$-saturating set of size $4q^{\prime }+4$ consisting of four pairwise
skew lines of $PG(3,q^{\prime })\subset PG(3,q).$
\end{corollary}

\begin{IEEEproof}
We use Theorems \ref{Th4_quattro} and \ref{Th4_dim6}.
\end{IEEEproof}

We now give an inductive construction of $v$-fold strong
blocking sets in $ PG(v,q)$. \medskip

\textbf{Construction A. }Let $H\cong PG(v,q)$ be a hyperplane
in $PG(v+1,q)$ , and let $B\subset H$ be a $v$-fold strong
blocking set in $H$. Let $ P_{1},P_{2},\ldots ,P_{v+1}$ be
$v+1$ independent points in $H$, and let $ l_1,\ldots,l_{v+1}$
be concurrent lines in $PG(v,q)$ such that $l_{i}\cap H=P_{i}$
for each $i$. Let
\begin{equation}
B^{\star }=B\cup \bigcup_{i=1,\ldots ,v+1}(l_{i}\setminus \{P_{i}\})\,.
\label{form4_inductive_strong}
\end{equation}

\begin{theorem}
\label{Th4_indu} Let $B$ be a $v$-fold strong blocking set in $PG(v,q)$ of
size$~k$. Then the set $B^{\star }$ of Construction A is a $(v+1)$-fold
strong blocking set in $PG(v+1,q)$ of size $k+1+(v+1)(q-1)$.
\end{theorem}

\begin{IEEEproof}
Let $H$ be the hyperplane in $PG(v+1,q)$ as in Construction A.
Let $H_{1}$ be any hyperplane in $PG(v+1,q)$. We need to show
that $H_{1}$ is generated by $v+1$ points in $B^{\star }$. When
$H=H_{1}$, this follows from the fact that $B$ must contain
$v+1$ independent points. Assume then that $H\neq H_{1}$, and
let $\Sigma =H\cap H_{1}$. As $\Sigma $ is a hyperplane in $H$,
there exist $v$ points $Q_{1},\ldots ,Q_{v}$ in $B$ which
generate $\Sigma $. Note that $\Sigma $ does not pass through a
point $ P_{i_{0}}$ for some $i_{0}\in \{1,\ldots ,v+1\}$, as
otherwise $\Sigma $ would coincide with $H$. Let $Q=H_{1}\cap
l_{i_{0}}$. As $Q\notin \Sigma $, and as $\Sigma $ is a
hyperplane of $H_{1}$, we have that $H_{1}=<\Sigma
,Q>=<Q_{1},\ldots ,Q_{v},Q>.$, with $\{Q_{1},\ldots
,Q_{v},Q\}\subset B^{\ast }$. This proves that $B^\star$ is a
$(v+1)-$fold blocking set. The size of $B^{\star }$ can be
easily calculated from~(\ref {form4_inductive_strong}).
\end{IEEEproof}

\begin{corollary}
\label{Cor4_Nfold_indu}In $PG(v,q)$, $v\geq 3$, there exists a $v$-fold
strong blocking set of size
\begin{equation}
(q-1)\left( \frac{v(v+1)}{2}-2\right) +v+5.
\end{equation}
\end{corollary}

\begin{IEEEproof}
By Theorem \ref{Th4_dim6}, in $PG(3,q)$ there exists a $3$-fold
strong blocking set of size $4q+4$. Then the assertion follows
by Theorem \ref {Th4_indu}, taking into account that
$4q+4+1+4(q-1)+1+5(q-1)+\ldots +1+v(q-1)=4q+4+(v-3)+(q-1)\left(
v(v+1)/2-6\right).$
\end{IEEEproof}

{}From Theorem \ref{Th4_quattro} we deduce the following result.

\begin{corollary}
\label{Cor4_PG(ro+1,q)_ro-sat}Let $q=(q^{\prime })^{\rho +1}$,
$\rho \geq 2$ . Then there exists a $\rho $-saturating set in
$PG(\rho +1,q)$ of size
\begin{equation}
(\sqrt[\rho +1]{q}-1)\left( \frac{(\rho +1)(\rho +2)}{2}-2\right) +\rho +6.
\label{form4_ro-sat}
\end{equation}
\end{corollary}

Note that the smallest previously known $\rho $-saturating sets
in $PG(\rho +1,$ $(q^{\prime })^{\rho +1}),$ $\rho \geq 2,$
have size $n=\frac{1}{2}(\sqrt[\rho +1]{q }-1)(\rho +1)(\rho
+2)+\rho +2$ \cite[Th.\thinspace 6]{DavO4}$,$ e.g.
$n=6q^{\prime }-2$ for $\rho =2$ and $n=10q^{\prime }-5$ for
$\rho =3$; from (\ref{form4_ro-sat}) we obtain sizes
$4q^{\prime }+4$ and $8q^{\prime }+1,$ respectively.

\begin{remark}
\label{Rem4_points-on-line} The codes associated to the
saturating sets of Corollaries \ref{Cor4_3fold_PG(3,q)}
and~\ref{Cor4_PG(ro+1,q)_ro-sat} will be used as starting codes
for $q^{m}$-concatenating constructions, see Sections
\ref{Sec7_R=3} and \ref{Sec8_R>=4}. Therefore, we need to treat
such codes as $(R,\ell)$-objects with $\ell >0$ and to obtain
the corresponding $(R,\ell )$-partitions, see Definition
\ref{Def2_R,l code}. To this end, it is useful to represent
some point $P_{i}$ of a line $l$ in $ PG(v,q)$ as a linear
combination with nonzero coefficients of $u$ other points
of~$l$. We compute some of the admissible values of $u$. Let
$l=\{P_{0},P_{1},\ldots ,P_{q}\}$. Without loss of generality
we identify $l$ with the projective line $PG(1,q)$, and assume
that $P_{0}=(0,1),$ $P_{1}=(1,0),$ $ P_{i}=(1,b_{i}),$ $i\geq
2,$ where $\{b_{2},\ldots ,b_{q}\}=F_{q}^{\ast }.$

\emph{i) }Clearly, for each $i=0,1,\ldots ,\left\lfloor
(q-2)/2\right\rfloor$ , the point $P_i$ can be written as
$P_{i}=c_{2i+1}P_{2i+1}+c_{2i+2}P_{2i+2} $, for some
$c_{2i+1},c_{2i+2}\in F_q^\ast$. So, $
P_{0}=c_{1}P_{1}+c_{2}P_{2}=c_{1}c_{3}P_{3}+c_{1}c_{4}P_{4}+c_{2}P_{2}=
c_{1}c_{3}P_{3}+c_{1}c_{4}P_{4}+c_{2}c_{5}P_{5}+c_{2}c_{6}P_{6},
$ and so on. Therefore, each $u\in \{2,3,\ldots ,\left\lfloor
(q+2)/2\right\rfloor\} $ is admissible.

\emph{ii)} Note that $P_1=(1,0)=-\sum_{i=2}^{q}(1,b_{i}). $ Then $u=q-1$ is
admissible.

\emph{iii) }Let $q\geq u\geq 3,$ $q\geq 4.$ Then, for any
$d_{i}\in F_q^\ast, $ one can always choose $a_{0}$ and $a_{1}$
in $F_q^\ast$ so that $
a_{0}(0,1)+a_{1}(1,0)+\sum_{i=2}^{u-1}d_{i}(1,b_{i})=a_{2}(1,b_{q})$
with some $a_{2}\in F_q^\ast.\smallskip $
\end{remark}

\subsection{ Small $\rho $-Saturating Sets in Spaces $
PG(v,(q^{\prime })^{\rho +1}),$ $v=\rho +2,\rho +3,\ldots
,2\rho -1$}

\begin{lemma}
\label{Lem4_key} Fix $1\leq k\leq v-1$. Let $B_{k}$ be the
subset of $ PG(v,q) $ consisting of the points whose Hamming
weight is at most $v-k+1,$ i.e. $B_{k}$ is the union of the
${\binom{v+1 }{k}}$ subspaces of equation $ x_{i_{1}}=\ldots
=x_{i_{k}}=0$, where $0\le i_1<i_2<\ldots<i_k\le v$. Then $
B_{k}$ is a $(k+1)$-strong blocking set.
\end{lemma}

\begin{IEEEproof}
Let $W$ be any $k$-dimensional subspace of $PG(v,q)$. Let $w_{1},\ldots
,w_{k+1}$ be $k+1$ independent points of $W$. Consider the matrix
\begin{equation*}
A_{W}=\left[ \renewcommand{\arraystretch}{0.7}
\begin{array}{c}
----w_{1}---- \\
----w_{2}---- \\
\vdots \\
----w_{k+1}----
\end{array}
\right]
\end{equation*}
whose rows are homogenous coordinates of points $w_{1},\ldots
,w_{k+1}$. As the rank of $A_{W}$ is equal to $k+1$, there
exists a non singular $ (k+1)\times (k+1)$ matrix $M=(m_{ij})$
such that $MA_{W}$ contains a submatrix $I_{k+1}$. Note that
the rows of $MA_{W}$ are the coordinates of $ (k+1)$ points of
$W$; more precisely the $i^{th}$-row of $MA_{W}$ is $
m_{i1}w_{1}+m_{i2}w_{2}+\ldots +m_{i(k+1)}w_{k+1}.$ Clearly
these points are independent, and they are contained in $B_{k}$
as $I_{k+1}$ is a submatrix of $MA_{W}$.
\end{IEEEproof}

Note that in the previous lemma
\begin{eqnarray}
|B_{k}|=\frac{1}{q-1}\sum_{i=1}^{v-k+1}(q-1)^{i}\binom{v+1}{i}=\notag\\
\frac{
V_{q}(v+1,v-k+1)-1}{q-1},  \label{form4_B_k}
\end{eqnarray}
see (\ref{form1_sphere}). Therefore the order of magnitude of
the size of $ B_{k}$ is $\binom{v+1}{k}q^{v-k}.$

\begin{theorem}
\label{Th4_Napoleon}Let $\rho $ be a positive integer. Let
$q=(q^{\prime })^{\rho +1}$ and $v>\rho +1$. Then in $PG(v,q)$
there exists a $\rho $ -saturating set of size
\begin{equation}
\frac{V_{q^{\prime }}(v+1,v-\rho +1)-1}{q^{\prime }-1}\sim \binom{v+1}{\rho }
q^{\frac{v-\rho }{\rho +1}}.  \label{form4_ro-sat_general}
\end{equation}
\end{theorem}

\begin{IEEEproof}
By Theorem \ref{Th4_quattro}, the set $B_{\rho }\subset
PG(v,q^{\prime })$, where $B_{\rho}$ is defined as in Lemma
\ref{Lem4_key}, is the desired $\rho $-saturating set.
\end{IEEEproof}

For some values of $v$ and $\rho $, the coefficient $\binom{v+1}{\rho }$ can
be improved. We show that this is possible for $v=4$, $\rho =2$.

Let $q=(q^{\prime })^{3}$. Let $E_{0}=(1,0,0,0,0)$,
$E_{1}=(0,1,0,0,0)$, $ E_{2}=(0,0,1,0,0)$, $E_{3}=(0,0,0,1,0)$,
and $E_{4}=(0,0,0,0,1)$ be points in $PG(4,q)$. For $k,i,j\in
\{0,1,2,3,4\}$, $k<i<j $, let $\pi _{k,i,j}$ be the plane in
$PG(4,q)$ generated by $E_{k},\,E_{i}$ and $E_{j}$. Let $\pi
_{1}=\pi _{0,1,2},$ $\pi _{2}=\pi _{0,3,4},$ $\pi _{3}=\pi
_{0,1,3},$ $\pi _{4}=\pi _{0,2,4},$ $\pi _{5}=\pi _{0,1,4},$
$\pi _{6}=\pi _{1,2,3},$ $\pi _{7}=\pi _{1,2,4},$ $\pi _{8}=\pi
_{1,3,4},$ $\pi _{9}=\pi _{2,3,4}.$ Let
\begin{equation}
S=\left( \bigcup_{s=1}^{9}\pi _{s}\right) \cap PG(4,q^{\prime }).
\label{form4_R=3_r=3t+2_2sat-9planes}
\end{equation}
The union $S$ of the nine planes $\pi_i$ consists of all points
of $ PG(4,q^{\prime }) $, apart from those belonging to the
following three disjoint classes: points with all non-zero
coordinates; points with precisely one zero coordinate; points
$(x,0,y,z,0)$ with $xyz\neq 0$. Therefore,
\begin{equation*}
|S|=\theta _{5,q^{\prime }}-(q^{\prime }-1)^{4}-5(q^{\prime
}-1)^{3}-(q^{\prime }-1)^{2}=9\sqrt[3]{q^{2}}-8\sqrt[3]{q}+4.
\end{equation*}

\begin{theorem}
\label{Th4_R=3_r=2t+2_9planes} Let $q=(q^{\prime })^{3}$\emph{.
}The set $S$ as in $(\ref{form4_R=3_r=3t+2_2sat-9planes})$ has
size $9\sqrt[3]{q^{2}}-8 \sqrt[3]{q}+4$, and it is a
$2$-saturating set in $PG(4,q)$.
\end{theorem}

\begin{IEEEproof}
Let $P$ be a point in $PG(4,q)$. Let $\pi $ be any plane of
$PG(4,q^{\prime })$ containing the subspace generated by
$P,P^{q^{\prime }},P^{(q^{\prime })^{2}}$. Clearly $\pi $
passes through $P$. Assume that $\pi $ does not pass through
$E_{0}$. Then among the points in $\{\pi \cap \pi _{s}\mid
s=1,\ldots ,5\}$ there are at least three non-collinear points
of $S$. Assume that $\pi $ passes through $E_{0}$. Let $H_{0}$
be the hyperplane generated by $E_{1},\ldots ,E_{4}$. Then $\pi
\cap H_{0}$ consists of a line $\ell $. Obviously, $\ell $
meets $\cup _{i=6}^{9}\pi _{i}$ in at least two non-collinear
points. Then $\pi $ passes through $3$ non-collinear points in
$S$.\smallskip
\end{IEEEproof}

\textbf{Open problem.} Reduce the coefficient $\binom{v+1}{\rho
}$ in (\ref{form4_ro-sat_general}), for generic $v$ and $\rho$.

\section{\label{Sec5_Tab}Tables of Upper Bounds on the Length Function
$\ell _{q}(r,R)$\emph{\ }for small $r$ and $R$}

We give tables of the values of$\;\overline{\ell }_{q}(r,R),$
i.e., the \emph{smallest known }lengths of a $q$-ary linear
code with codimension $r$ and covering radius $R.$ Obviously,
$\ell _{q}(r,R)\leq \overline{\ell } _{q}(r,R)$ holds. The dot
\textquotedblleft $\centerdot $\textquotedblright\ appears in a
table when $\ell _{q}(r,R)=\overline{\ell }_{q}(r,R)$ holds.
Subscripts indicate the minimum distance $d$ of the
$[\,\overline{\ell } _{q}(r,R),\overline{\ell
}_{q}(r,R)-r,d]_{q}R$ codes. Multiple subscripts mean that the
value of $\overline{\ell }_{q}(r,R)$ is provided by codes with
distinct distances.

Table \ref{table1} gives values of $\overline{\ell }_{q}(3,2).$
We used \cite[ Tabs\thinspace
2,4]{DMP-JCTA},\cite[Tab.\thinspace I]{DMPIEEE},\cite[
Tab.\thinspace 3]{DFMP-JG}, Theorem \ref{Th4_q4-1sat},
Corollary \ref{Cor4_q6-1sat}, the relation $\ell _{(q^{\prime
})^{2}}(3,2)\leq 3q^{\prime }-1$ \cite[Th.\thinspace
5.2]{Dav95}, and computer search made in this work. Note that
the distance $d=4$ occurs when the code arises from a complete
arc in the plane $PG(2,q)$.

\begin{table*}
\caption{Upper Bounds $\overline{\ell }_{q}=\overline{\ell
}_{q}(3,2)$ on the Length Function $\ell _{q}(3,2)$}
\renewcommand{\arraystretch}{1.0}
\begin{equation*}
\begin{array}
{|@{\,\,}r@{\,\,\,}r@{\,\,}|@{\,\,}r@{\,\,\,}r@{\,\,}|@{\,\,}r@{\,\,\,}r@{\,\,}|@{\,\,}r@{\,\,\,}r@{\,\,}|@{\,\,}r@{\,\,\,}
r@{\,\,}|@{\,\,}r@{\,\,\,}r@{\,\,}|@{\,\,}r@{\,\,\,}r@{\,\,}|@{\,\,}r@{\,\,\,}r@{\,\,}|@{\,\,}r@{\,\,\,}r@{\,\,}|}
\hline
q & \overline{\ell }_{q}\phantom{\overline{\overline{H}}} & q & \overline{\ell }_{q}
\phantom{_{3,4}} & q & \overline{\ell }_{q}\phantom{_{3,4}} & q & \overline{
\ell }_{q}\phantom{_{3}} & q & \overline{\ell }_{q}\phantom{_{3}} & q &
\overline{\ell }_{q}\phantom{_{3}} & q & \overline{\ell }_{q}\phantom{_{3}}
& q & \overline{\ell }_{q}\phantom{_{3}} & q & \overline{\ell }_{q}
\phantom{_{3}} \\ \hline
3 & 4_{4\phantom{,4}}\centerdot & 64 & 19_{3\phantom{,4}} & 167 & 42_{3,4} &
283 & 58_{3} & 431 & 75_{3} & 577 & 89_{3} & 729 & 80_{3} & 887 & 116_{3} &
1051 & 128_{3} \\
4 & 5_{3\phantom{,4}}\centerdot & 67 & 23_{3,4} & 169 & 38_{3\phantom{,4}} &
289 & 50_{3} & 433 & 75_{3} & 587 & 90_{3} & 733 & 102_{3} & 907 & 117_{3} &
1061 & 129_{3} \\
5 & 6_{3,4}\centerdot & 71 & 22_{4\phantom{,4}} & 173 & 42_{3\phantom{,4}} &
293 & 59_{3} & 439 & 75_{3} & 593 & 90_{3} & 739 & 103_{3} & 911 & 118_{3} &
1063 & 129_{3} \\
7 & 6_{3,4}\centerdot & 73 & 24_{4\phantom{,4}} & 179 & 43_{3\phantom{,4}} &
307 & 60_{3} & 443 & 76_{3} & 599 & 91_{3} & 743 & 104_{3} & 919 & 118_{3} &
1069 & 129_{3} \\
8 & 6_{4\phantom{,4}}\centerdot & 79 & 26_{3,4} & 181 & 43_{3\phantom{,4}} &
311 & 61_{3} & 449 & 76_{3} & 601 & 91_{3} & 751 & 105_{3} & 929 & 119_{3} &
1087 & 130_{3} \\
9 & 6_{4\phantom{,4}}\centerdot & 81 & 26_{3,4} & 191 & 45_{3\phantom{,4}} &
313 & 61_{3} & 457 & 77_{3} & 607 & 91_{3} & 757 & 105_{3} & 937 & 120_{3} &
1091 & 131_{3} \\
11 & 7_{4\phantom{,4}}\centerdot & 83 & 26_{3\phantom{,4}} & 193 & 45_{3
\phantom{,4}} & 317 & 62_{3} & 461 & 77_{3} & 613 & 92_{3} & 761 & 105_{3} &
941 & 120_{3} & 1093 & 131_{3} \\
13 & 8_{4\phantom{,4}}\centerdot & 89 & 28_{3,4} & 197 & 46_{3\phantom{,4}}
& 331 & 63_{3} & 463 & 77_{3} & 617 & 92_{3} & 769 & 106_{3} & 947 & 121_{3}
& 1097 & 131_{3} \\
16 & 9_{3,4}\centerdot & 97 & 29_{3\phantom{,4}} & 199 & 46_{3\phantom{,4}}
& 337 & 64_{3} & 467 & 78_{3} & 619 & 92_{3} & 773 & 106_{3} & 953 & 121_{3}
& 1103 & 131_{3} \\
17 & 10_{3,4}\phantom{\centerdot} & 101 & 30_{3,4} & 211 & 48_{3\phantom{,4}}
& 343 & 64_{3} & 479 & 79_{3} & 625 & 62_{3} & 787 & 107_{3} & 961 & 92_{3}
& 1109 & 132_{3} \\
19 & 10_{3,4}\phantom{\centerdot} & 103 & 30_{3\phantom{,4}} & 223 & 49_{3
\phantom{,4}} & 347 & 65_{3} & 487 & 80_{3} & 631 & 94_{3} & 797 & 108_{3} &
967 & 123_{3} & 1117 & 132_{3} \\
23 & 10_{4\phantom{,4}}\centerdot & 107 & 31_{3\phantom{,4}} & 227 & 50_{3
\phantom{,4}} & 349 & 65_{3} & 491 & 81_{3} & 641 & 95_{3} & 809 & 109_{3} &
971 & 123_{3} & 1123 & 133_{3} \\
25 & 12_{3,4}\phantom{\centerdot} & 109 & 31_{3\phantom{,4}} & 229 & 50_{3
\phantom{,4}} & 353 & 66_{3} & 499 & 81_{3} & 643 & 95_{3} & 811 & 110_{3} &
977 & 124_{3} & 1129 & 133_{3} \\
27 & 12_{3,4}\phantom{\centerdot} & 113 & 32_{3\phantom{,4}} & 233 & 51_{3
\phantom{,4}} & 359 & 66_{3} & 503 & 82_{3} & 647 & 95_{3} & 821 & 110_{3} &
983 & 124_{3} & 1151 & 135_{3} \\
29 & 13_{3,4}\phantom{\centerdot} & 121 & 32_{3\phantom{,4}} & 239 & 51_{3
\phantom{,4}} & 361 & 56_{3} & 509 & 82_{3} & 653 & 96_{3} & 823 & 110_{3} &
991 & 124_{3} & 1153 & 135_{3} \\
31 & 14_{3,4}\phantom{\centerdot} & 125 & 34_{3\phantom{,4}} & 241 & 52_{3
\phantom{,4}} & 367 & 67_{3} & 512 & 82_{3} & 659 & 96_{3} & 827 & 110_{3} &
997 & 125_{3} & 1163 & 135_{3} \\
32 & 13_{3}\phantom{_{,4}\centerdot} & 127 & 35_{3,4} & 243 & 52_{3
\phantom{,4}} & 373 & 68_{3} & 521 & 84_{3} & 661 & 96_{3} & 829 & 110_{3} &
1009 & 124_{3} & 1171 & 136_{3} \\
37 & 15_{4}\phantom{_{,3}\centerdot} & 128 & 34_{3,4} & 251 & 53_{3
\phantom{,4}} & 379 & 69_{3} & 523 & 83_{3} & 673 & 98_{3} & 839 & 111_{3} &
1013 & 125_{3} & 1181 & 137_{3} \\
41 & 16_{4}\phantom{_{,3}\centerdot} & 131 & 35_{3\phantom{,4}} & 256 & 42_{3
\phantom{,4}} & 383 & 69_{3} & 529 & 68_{3} & 677 & 98_{3} & 853 & 113_{3} &
1019 & 126_{3} & 1187 & 137_{3} \\
43 & 16_{4}\phantom{_{,3}\centerdot} & 137 & 36_{3\phantom{,4}} & 257 & 54_{3
\phantom{,4}} & 389 & 70_{3} & 541 & 85_{3} & 683 & 99_{3} & 857 & 113_{3} &
1021 & 126_{3} & 1193 & 138_{3} \\
47 & 18_{3,4}\phantom{\centerdot} & 139 & 37_{3,4} & 263 & 55_{3\phantom{,4}}
& 397 & 71_{3} & 547 & 86_{3} & 691 & 99_{3} & 859 & 113_{3} & 1024 & 95_{3}
& 1201 & 138_{3} \\
49 & 18_{4}\phantom{_{,3}\centerdot} & 149 & 39_{3,4} & 269 & 56_{3
\phantom{,4}} & 401 & 71_{3} & 557 & 87_{3} & 701 & 100_{3} & 863 & 114_{3}
& 1031 & 127_{3} & 1213 & 139_{3} \\
53 & 18_{4}\phantom{_{,3}\centerdot} & 151 & 39_{3,4} & 271 & 56_{3
\phantom{,4}} & 409 & 72_{3} & 563 & 87_{3} & 709 & 101_{3} & 877 & 115_{3}
& 1033 & 127_{3} & 1217 & 139_{3} \\
59 & 20_{4}\phantom{_{,3}\centerdot} & 157 & 40_{3,4} & 277 & 57_{3
\phantom{,4}} & 419 & 73_{3} & 569 & 88_{3} & 719 & 102_{3} & 881 & 115_{3}
& 1039 & 127_{3} &  &  \\
61 & 20_{4}\phantom{_{,3}\centerdot} & 163 & 41_{3,4} & 281 & 57_{3
\phantom{,4}} & 421 & 73_{3} & 571 & 88_{3} & 727 & 102_{3} & 883 & 115_{3}
& 1049 & 128_{3} &  &  \\ \hline
\end{array}
\end{equation*}
\label{table1}
\end{table*}
From Table \ref{table1} the following result is obtained.

\begin{theorem}
\label{Th5_1sat_byTable}For the length function $\ell _{q}(3,2)$,
\begin{eqnarray}
\ell _{q}(3,2)\leq a_{q}\sqrt{q},&\text{~~with }a_{q}<3\text{ if }q\leq 109,\notag\\
\text{ }a_{q}<3.5\text{ if }q\leq 349,&\text{ }a_{q}<4\text{ if }q\leq 1217.
\label{form5_1satplane}
\end{eqnarray}
\end{theorem}

In Table \ref{table2} we give a number of concrete sizes of
1-saturating sets and complete caps in $PG(2,q),$ $q=p^{2t+1},$
taken from \cite[Tab.\thinspace 2]
{Giul-plane},\cite[Ap.,\thinspace Lem.\thinspace 4.3]{Gi2}, and
\cite[ Tab.\thinspace 1]{DGMP-Submit}. These sizes are the
values of $\overline{ \ell }_{q}(3,2).$

\begin{table*}
\caption{Upper Bounds $\overline{\ell }_{q}=\overline{\ell
}_{q}(3,2)$ on the Length Function $\ell _{q}(3,2)$ for
$q=p^{2t+1}$}
\renewcommand{\arraystretch}{1.1}
\begin{equation*}
\begin{array}{|@{\,\,}r@{\,\,}r@{\,\,}|@{\,\,}l@{\,\,}r@{\,\,}|@{\,\,}l@{\,\,}r@{\,\,}|@{\,\,}r@{\,\,}r@{\,\,}
|@{\,\,}r@{\,\,}r@{\,\,}|@{\,\,}r@{\,\,}r@{\,\,}|@{\,\,}r@{\,\,}r@{\,\,}|@{\,\,}rr@{\,\,}|}
\hline
q & \overline{\ell }_{q} & q & \overline{\ell }_{q} & q & \overline{\ell }
_{q} & q & \overline{\ell }_{q} & q & \overline{\ell }_{q} & q & \overline{
\ell }_{q} & q & \overline{\ell }_{q} & q \phantom{\overline{\overline{H}}}& \overline{\ell }_{q} \\ \hline
2^{11} & 201_{4} & 2^{17} & 2576_{3} & 3^{9} & 764_{3} & 5^{5} & 376_{3} &
7^{5} & 1030_{3} & 11^{5} & 3994_{3} & 13^{7} & 85712_{3} & 19^{5}\phantom{\overline{\overline{H}}} &
20578_{3} \\
2^{13} & 461_{4} & 2^{19} & 5210_{3} & 3^{11} & 2771_{3} & 5^{7} & 1877_{3}
& 7^{7} & 7205_{3} & 11^{7} & 43947_{3} & 17^{5} & 14740_{3} &  &  \\
2^{15} & 993_{4} & 3^{7} & 245_{3} & 3^{13} & 8788_{3} & 5^{9} & 9609_{3} &
7^{9} & 50947_{3} & 13^{5} & 6592_{3} & 17^{7} & 250599_{3} &  &  \\ \hline
\end{array}
\end{equation*}
\label{table2}
\end{table*}

Using \cite[Tab.\thinspace 1]{DavO4},\cite[Tabs
II,III]{DMPIEEE},\cite[ Tabs\thinspace III-V]{DFMP-IEEE-LO},
Theorem \ref{Th4_dim6} and Corollary
\ref{Cor4_PG(ro+1,q)_ro-sat}, we obtained Table \ref{table3}
where values of $ \overline{\ell }_{q}(4,3)$ are listed. The
distances $d=4$ and $d=5$ occur, respectively, when the code
arises from an incomplete cap and a complete arc in $PG(3,q)$
\cite{DMPIEEE},\cite{DFMP-IEEE-LO}.

\begin{table*}
\caption{Upper Bounds $\overline{\ell
}_{q}=\overline{\ell }_{q}(4,3)$ on the Length Function $\ell
_{q}(4,3)$}
\renewcommand{\arraystretch}{1.0}
\begin{equation*}
\begin{array}{|@{\,\,}r@{\,\,\,}l@{\,\,}|@{\,\,}r@{\,\,\,}l@{\,\,}|@{\,\,}r@{\,\,\,}l@{\,\,}|@{\,\,}
r@{\,\,\,}l@{\,\,}|@{\,\,}r@{\,\,\,}l@{\,\,}|@{\,\,}r@{\,\,\,}l@{\,\,}|@{\,\,}@{\,}r@{\,\,\,}l@{\,\,}
|@{\,\,}r@{\,\,\,}l@{\,\,}|@{\,\,}r@{\,\,\,}l@{\,\,}|}
\hline
q & \phantom{1}\overline{\ell }_{q}\phantom{\overline{\overline{H}}} & q & \phantom{1}\overline{\ell }_{q} & q
& \phantom{1}\overline{\ell }_{q} & q & \phantom{1}\overline{\ell }_{q} & q
& \phantom{1}\overline{\ell }_{q} & q & \phantom{1}\overline{\ell }_{q} & q
& \phantom{1}\overline{\ell }_{q} & q & \phantom{1}\overline{\ell }_{q} & q
& \phantom{1}\overline{\ell }_{q} \\ \hline
2 & \phantom{1}5_{3,4}\centerdot & 27 & 11_{3,4,5} & 71 & 16_{4,5} & 127 &
21_{3,5} & 191 & 25_{3,5} & 257 & 28_{3,5} & 337 & 31_{3} & 409 & 34_{3,5} &
491 & 36_{4} \\
3 & \phantom{1}5_{4,5}\centerdot & 29 & 11_{3,4,5} & 73 & 16_{4} & 128 &
21_{3,5} & 193 & 25_{3,5} & 263 & 28_{3,5} & 343 & 31_{4} & 419 & 34_{3} &
499 & 37_{3,5} \\
4 & \phantom{1}5_{5}\centerdot & 31 & 11_{4} & 79 & 17_{3,5} & 131 & 21_{3,5}
& 197 & 25_{3,5} & 269 & 29_{3,5} & 347 & 32_{3,5} & 421 & 34_{3} & 503 &
37_{3,5} \\
5 & \phantom{1}6_{3,4,5}\centerdot & 32 & 12_{3,4,5} & 81 & 17_{4} & 137 &
22_{3,5} & 199 & 25_{5} & 271 & 29_{3,5} & 349 & 32_{3,5} & 431 & 35_{3,5} &
509 & 37_{3,5} \\
7 & \phantom{1}7_{3,4}\centerdot & 37 & 12_{4,5} & 83 & 17_{4} & 139 &
22_{3,5} & 211 & 26_{3,5} & 277 & 29_{3,5} & 353 & 32_{3,5} & 433 & 35_{3} &
8^{3} & 36_{3} \\
8 & \phantom{1}7_{3,4,5}\centerdot & 41 & 13_{3,4,5} & 89 & 18_{3,5} & 149 &
22_{3,5} & 223 & 27_{3,5} & 281 & 29_{3,5} & 359 & 32_{3,5} & 439 & 35_{3,5}
& 521 & 37_{4} \\
9 & \phantom{1}7_{4}\centerdot & 43 & 13_{4,5} & 97 & 19_{3,5} & 151 & 22_{4}
& 227 & 27_{3,5} & 283 & 29_{3,5} & 361 & 32_{3} & 443 & 35_{3,5} & 523 &
38_{3} \\
11 & \phantom{1}8_{3,4,5}\centerdot & 47 & 14_{3,4,5} & 101 & 19_{3.5} & 157
& 23_{3,5} & 229 & 27_{3,5} & 289 & 29_{4} & 367 & 32_{4} & 449 & 35_{3,5} &
529 & 38_{5} \\
13 & \phantom{1}8_{4,5} & 49 & 14_{3,4,5} & 103 & 19_{3,5} & 163 & 23_{5} &
233 & 27_{3,5} & 293 & 29_{4} & 373 & 33_{3,5} & 457 & 35_{4} & 541 & 38_{5}
\\
16 & \phantom{1}9_{3,4,5} & 53 & 15_{3,4,5} & 107 & 19_{4} & 167 & 24_{3,5}
& 239 & 27_{3,5} & 307 & 30_{3,5} & 379 & 33_{3,5} & 461 & 36_{3,5} & 547 &
38_{4} \\
17 & \phantom{1}9_{3,4,5} & 59 & 15_{3,4,5} & 109 & 20_{3,5} & 169 & 24_{3,5}
& 241 & 28_{3,5} & 311 & 30_{4} & 383 & 33_{3,5} & 463 & 36_{3} & 557 &
39_{5} \\
19 & \phantom{1}9_{4,5} & 61 & 15_{4} & 113 & 20_{3,5} & 173 & 24_{3,5} & 243
& 28_{3,5} & 313 & 30_{4} & 389 & 33_{4} & 467 & 36_{3} & 563 & 39_{5} \\
23 & 10_{3,4,5} & 64 & 16_{3,4,5} & 121 & 20_{4} & 179 & 24_{3,5} & 251 &
28_{3,5} & 317 & 30_{4} & 397 & 34_{3,5} & 479 & 36_{3} & 9^{3} & 40_{3} \\
25 & 11_{3,4,5} & 67 & 16_{3,4,5} & 125 & 21_{3,5} & 181 & 24_{4} & 256 &
28_{3,5} & 331 & 31_{3,5} & 401 & 34_{3,5} & 487 & 36_{3,5} & 11^{3} & 48_{3}
\\ \hline
\end{array}
\end{equation*}
\label{table3}
\end{table*}

From Table \ref{table3} we obtain the following theorem.

\begin{theorem}
\label{Th5_2sat_r0=4_byTabIII}For the length function $\ell
_{q}(4,3)$,
\begin{eqnarray}
\ell _{q}(4,3)\leq b_{q}\sqrt[3]{q},&\;\text{with}\;b_{q}<4\text{ if }q\leq
83,\notag\\
b_{q}<4.5\text{ if }q\leq 343,&\;b_{q}<5\text{ if }q\leq 563.
\label{form5_befend}
\end{eqnarray}
\end{theorem}

In Table \ref{table4} the values of $\overline{\ell }_{q}(5,3)$
are given. We use \cite
[Tab.\thinspace1]{DavO4},\cite[Tabs\thinspace
III,IV]{DFMP-IEEE-LO} for $ q\leq 7$ and the computer search
made in this work for $8\leq q\leq 32.$ For $37\leq q\leq 43,$
we apply the direct sum (see Section \ref{Sec6_R=2}) of the
$[\overline{\ell }_{q}(3,2),\overline{\ell }_{q}(3,2)-3]_{q}2$
code of Table I and the $[q+1,q-1]_{q}1$ Hamming code. The
distances $d=4$ and $d=6$ occur note, respectively, when the
code arises from an incomplete cap and an arc in $PG(4,q)$.

\begin{table*}
\caption{Upper Bounds $\overline{\ell }_{q}=\overline{\ell
}_{q}(5,3)$ on the Length Function $\ell _{q}(5,3)$}
\renewcommand{\arraystretch}{1.0}
\begin{equation*}
\renewcommand{\arraystretch}{1.1}
\begin{array}{|@{\,\,}r@{\,\,\,}l@{\,\,}|@{\,\,}r@{\,\,\,}l@{\,\,}|@{\,\,}
r@{\,\,\,}l|@{\,\,}r@{\,\,\,}l@{\,\,}|@{\,\,}r@{\,\,\,}l@{\,\,}|}
\hline
q & \overline{\ell }_{q}\phantom{\overline{\overline{H}}} & q & \overline{\ell }_{q} & q & \overline{\ell }
_{q} & q & \overline{\ell }_{q} & q & \overline{\ell }_{q} \\ \hline
2 & 6_{5,6}\centerdot & 8 & 14_{3,4} & 17 & 25_{3,4} & 29 & 38_{3,4} & 43 &
60_{3} \\
3 & 8_{3,4}\centerdot & 9 & 16_{3,4} & 19 & 27_{4} & 31 & 40_{3,4} &  &  \\
4 & 9_{3,4}\centerdot & 11 & 18_{3,4} & 23 & 32_{3,4} & 32 & 41_{3,4} &  &
\\
5 & 10_{4}\centerdot & 13 & 21_{3,4} & 25 & 34_{3,4} & 37 & 54_{3} &  &  \\
7 & 13_{3,4} & 16 & 24_{3,4} & 27 & 36_{3,4} & 41 & 58_{3} &  &  \\ \hline
\end{array}
\end{equation*}
\label{table4}
\end{table*}

From Table \ref{table4} we obtain the following theorem.

\begin{theorem}
\label{Th5_2sat_r0=5_byTabIV}For the length function $\ell
_{q}(5,3)$,
\begin{eqnarray*}
\ell _{q}(5,3)\leq c_{q}\sqrt[3]{q^{2}},&\;\text{with}\;c_{q}<4\text{ if }
q\leq 27,\;\\
c_{q}<4.2\text{ if }q\leq 32,&\;c_{q}<5\text{ if }q\leq 43.
\end{eqnarray*}
\end{theorem}

\section{\label{Sec6_R=2}Codes with Covering radius $R=2$}

\subsection{ Direct sum and doubling constructions}

The \emph{direct sum }construction (DS) forms an $
[n_{1}+n_{2},n_{1}+n_{2}-(r_{1}+r_{2})]_{q}R$ code $V$ with
$R=R_{1}+R_{2}$ from two codes: an
$[n_{1},n_{1}-r_{1}]_{q}R_{1}$ code $V_{1}$ and an $
[n_{2},n_{2}-r_{2}]_{q}R_{2}$ code $V_{2}$
\cite{GrSl},\cite{Coh},\cite {Handbook-coverings}. The
parity-check matrix $\mathbf{H}$ of the new code $ V $ has the
form
\begin{equation*}
\mathbf{H=}\left[
\begin{array}{@{\,\,}c@{\,\,}c@{\,\,}}
\mathbf{H}_{1} & \mathbf{0}_{r_{1}} \\
\mathbf{0}_{r_{2}} & \mathbf{H}_{2}
\end{array}
\right]
\end{equation*}
where $\mathbf{H}_{1}$ and $\mathbf{H}_{2}$ are parity-check matrices of the
starting codes $V_{1}$ and $V_{2},$ respectively. Construction DS is denoted
by $\oplus ,$ i.e. $V_{1}\oplus V_{2}=V$ or
\begin{eqnarray}
\lbrack n_{1},n_{1}-r_{1}]_{q}R_{1}\oplus \lbrack
n_{2},n_{2}-r_{2}]_{q}R_{2}=\notag\\
\,[n_{1}+n_{2},n_{1}+n_{2}-(r_{1}+r_{2})]_{q}(R_{1}+R_{2}).
\label{form6_DS}
\end{eqnarray}
DS construction yields that
\begin{equation*}
\ell _{q}(r_{1}+r_{2},R_{1}+R_{2})\leq \ell _{q}(r_{1},R_{1})+\ell
_{q}(r_{2},R_{2}).
\end{equation*}

In \cite{davShumen} Construction CP1 (\textquotedblleft
codimension plus one\textquotedblright ) is proposed. The
construction is similar to the construction
in~\cite{Hon-length}. From an $[n,n-r]_{q}2$ code $V_{1}$
Construction CP1 forms an $[f_{q}(n),f_{q}(n)-(r+1)]_{q}2$ code
$V$ where $ f_{3}(n)=2n,$ $f_{4}(n)=3n-1,$ $f_{5}(n)=3n.$ For
$q=3$ Construction CP1 is a \emph{doubling construction.} In
this case the parity-check matrix $ \mathbf{H}$ of the new code
$V$ has the form
\begin{equation}
\mathbf{H=}\left[
\begin{array}{@{\,\,}c@{\,\,}c@{\,\,}}
\mathbf{0} & \mathbf{1} \\
\mathbf{H}_{1} & \mathbf{H}_{1}
\end{array}
\right] ,\text{ }q=3,  \label{form6_doubling}
\end{equation}
where $\mathbf{0}$ and $\mathbf{1}$ is the row of all zeroes
and units, respectively, and $\mathbf{H}_{1}$ is a parity-check
matrix of the starting code $V_{1}$. By~(\ref{form6_doubling}),
see also \cite{ostArs},
\begin{equation}
\ell _{3}(r+1,2)\leq 2\ell _{3}(r,2).  \label{form6_l3(r+1,2)}
\end{equation}

\subsection{Infinite Code Families of Even Codimension $r=2t$}

Let $q=3.$ By applying the doubling construction of
(\ref{form6_doubling}) to the codes of \cite[Th.\thinspace
1]{coc},\cite[Th.\thinspace 4]{DavO1} and by using the codes of
\cite[Th.\thinspace 11]{DavO1} we obtain an infinite family of
$[n,n-r]_{3}2$ codes with the following parameters
\begin{eqnarray}
\mathcal{A}_{2,3}^{(0)} &:&R=2,\text{ }r=2t\geq 4,\text{ }q=3,\text{ }r\neq
8,\text{ }\overline{\mu }_{3}(2)\approx \frac{25}{18},  \notag \\
n &=&\frac{5}{2}\cdot 3^{\frac{r-2}{2}}-\frac{1}{2}+\left\{
\begin{array}{l@{\,\,}l}
0 & \text{if }r=4c+2 \\
\frac{1}{2}\cdot 3^{\frac{r}{4}}-\frac{1}{2}\smallskip & \text{if }r=8c+4 \\
\frac{1}{2}\cdot 3^{\frac{r+4}{4}}-\frac{1}{2} & \text{if }r=8c
\end{array}
\right. .  \label{form6_R=2_r=2t_1}
\end{eqnarray}
For $r=4,$ from (\ref{form6_R=2_r=2t_1}) we obtain an $[8,4]_{3}2$ code.
Note that by \cite[Tab.\thinspace II]{BV}, $\ell _{3}(4,2)=8$ holds.

Let $q\geq 4.$ The geometrical constructions (named
\textquotedblleft oval plus line\textquotedblright ) give
$[2q+1,2q-3]_{q}2$ codes, see \cite[ p.\thinspace 104]{BrPlWi}
for even$~q$ and \cite[Th.\thinspace
3.1]{DavParis},\cite[Th.\thinspace 5.1]{Dav95} for arbitrary
$q$. By computer, using the back-tracking algorithms
\cite{MPAustr},\cite{DFMP-JG}, we have proved the following
proposition.

\begin{proposition}
\label{prop6_l4(4,2)} $\ell _{4}(4,2)=9$.
\end{proposition}

No examples of $[n,n-4]_{q}2$ codes with $n<2q+1,$ seem to be
known.

\textbf{Open problem.} To prove that $\ell _{q}(4,2)=2q+1$ for
$q\geq 5.$\smallskip

In \cite{DavO5} the parity-check matrices of the codes of
\cite{BrPlWi},\cite {Dav95} are modified and used as starting
$(R,\ell)$-objects in $q^{m}$-concatenating constructions. As a
result, an infinite family of $ [n,n-r]_{q}2$ codes is obtained
with the following parameters \cite[ Th.\thinspace 9]{DavO5}:
\begin{eqnarray}
\mathcal{A}_{2,q}^{(0)}&:&R=2,\text{ }r=2t\geq 4,\text{ }
q\geq 7,\text{ }q\neq 9,\text{
}r\neq 8,12,
\label{form6_R=2_r=2t_2}\\
\text{ }n&\,=&\,2q^{\frac{r-2}{2}}+q^{\frac{r-4}{2}},\text{ }
\overline{\mu }_{q}(2) <2-\frac{2}{q}-\frac{3}{2q^{2}}+\frac{1}{q^{3}}+
\frac{1}{q^{4}}. \notag
\end{eqnarray}
Also, in \cite{DavO5} codes with $r=8,12,$
$n=2q^{\frac{r-2}{2}}+q^{\frac{r-4
}{2}}+q^{\frac{r-6}{2}}+q^{\frac{r-8}{2}},$ $q\geq 7,$ $q\neq
9,$ are given.

For $q=4,5,9$ in \cite[Ex. 5]{DavNBCR2} an infinite family of $[n,n-r]_{q}2$
codes is obtained with
\begin{eqnarray}
\mathcal{A}_{2,q}^{(0)}&:&R=2,\text{ }r=2t\geq 4,\text{ }
q=4,5,9,\label{form6_R=2_r=2t_3}\\
n&\,=&\,2q^{\frac{r-2}{2}}+q^{\frac{r-4}{2}}+\left\lfloor q^{\frac{
r-6}{2}}\right\rfloor ,\text{ }r \neq 8,12,14,20\text{ if }q=4, \notag \\
r&\neq& 8,12\text{ if }q=5,9,\text{ }\overline{\mu }_{q}(2)<2-\frac{2}{q
}+\frac{1}{2q^{2}}-\frac{2}{q^{3}}+\frac{2}{q^{4}}.  \notag
\end{eqnarray}
Also, codes with $q=4,5,9,$ $n=2q^{3}+q^{2}+2q+2,$ $r=8,$ and $
n=q^{5}+\theta _{6,q},$ $r=12,$ are given.

\subsection{More on 1-Saturating Sets in Projective
Planes $ PG(2,q)$}

We recall here some of the known results on small
$1$-saturating sets in $ PG(2,q)$. (For the new 1-saturating
sets obtained in this paper we refer to Section
\ref{Sec4_NewSmSat} and Tables I,II of Section~\ref{Sec5_Tab}).

For large $q$ the existence of 1-saturating sets in $PG(2,q)$ of size at
most $5\sqrt{q\log q}$ was shown by means of \emph{probabilistic methods} in
\cite{Bo-Sz-Ti},\cite{kovS92}.

The following results are given by explicit constructions.

In $PG(2,q),$ $q=(q^{\prime })^{2},$ a\ 1-saturating set of
size $3\sqrt{q} -1 $ is obtained in \cite[Th.\ 5.2]{Dav95}.

In the plane $PG(2,q),$ $q=(q^{\prime })^{m},$ $m\geq 2,$
projectively non-equivalent 1-saturating sets of size
$2q^{\frac{m-1}{m}}+\sqrt[m]{q}$ are obtained in
\cite[Th.\thinspace 2]{DavO4},\cite[Th.\thinspace 3.2]
{Giul-plane}.

In \cite{Bartocci},\cite{Sz-survey89},\cite{GiuTor-Ars04},\cite{Bo-Sz-Ti}
1-saturating sets in $PG(2,q)$ of size approximately $cq^{\frac{3}{4}}$ with
a constant $c$ independent of $q$ are constructed.

In \cite{Giul-plane} constructions of 1-saturating $n$-sets in $PG(2,q)$ of
size $n$ about $3q^{\frac{2}{3}}$ are proposed. In particular the following
upper bounds on $n$ are obtained for $p$ prime:
\begin{equation}
n\leq \left\{ \renewcommand{\arraystretch}{1.4}
\begin{array}{l}
\frac{2q}{p^{t}}+\frac{(p^{t}-1)^{2}}{p-1}+1,~q=p^{m},\text{ }m\geq 2t; \\
\frac{2}{p}\sqrt[3]{(qp)^{2}}+\frac{\sqrt[3]{(qp)^{2}}-2\sqrt[3]{qp}+1}{p-1}
+1, ~ q=p^{3t-1}; \\
\min\limits_{v=1,\ldots ,2t+1}\left\{ (v+1)p^{t+1}+\frac{(p^{t}-1)^{2v}}{
(p-1)^{v}(p^{2t+1}-1)^{(v-1)}}+2\right\} ,\\  q=p^{2t+1}.
\end{array}
\right.   \label{form6_R=2_r=2t+1_q^2/3}
\end{equation}
Several triples $(t,p,v)$ such that $n<5\sqrt{q\log q}$ are
obtained in \cite {Giul-plane}.

\subsection{ Infinite Code Families of Odd Codimension $
r=2t+1$}

Let $q=3.$ By \cite[Th.\thinspace 1]{coc},\cite[Ths 4 and
9]{DavO1}, there exists an infinite family of $[n,n-r]_{3}2$
codes with the following parameters:
\begin{eqnarray}
\mathcal{A}_{2,3}^{(1)} &:&R=2,\text{ }r=2t+1\geq 5,\text{ }q=3,\text{ }
r\neq 7,\text{ }\overline{\mu }_{3}(2)\approx \frac{25}{24},  \notag \\
n &=&\frac{5}{4}\sqrt{3}\cdot 3^{\frac{r-2}{2}}-\frac{1}{4}+\left\{
\begin{array}{l@{\,\,}l}
0 & \text{if }r=4c+1 \\
\frac{3}{4}\cdot 3^{\frac{r+1}{4}}-\frac{3}{4}\smallskip & \text{if }r=8c+3
\\
\frac{3}{4}\cdot 3^{\frac{r+5}{4}}-\frac{3}{4} & \text{if }r=8c+7
\end{array}
\right. .  \label{form6_R=2_r=2t+1_q=3}
\end{eqnarray}

Let $q = 4$. In \cite{DavO3} an infinite family of
$[n,n-r]_{4}2$ codes is obtained with parameters
\begin{eqnarray}
\mathcal{A}_{2,4}^{(1)} &:&R=2,\text{ }r=2t+1\geq 5,\text{ }q=4,\text{ }
r\neq 7,11,13,19, \notag
\\
n &\,=\,&2\cdot 4^{\frac{r-2}{2}}+\frac{3}{2}\cdot 4^{\frac{r-4}{2}},
\text{ }\overline{\mu }_{4}(2)\approx 1.587 .
\label{form6_R=2_r=2t+1_q=4}
\end{eqnarray}

Let $q$ = 5. In \cite[Ths 5,10]{DavO1} an infinite family of
$[n,n-r]_{5}2$ codes is obtained with
\begin{eqnarray}
\mathcal{A}_{2,5}^{(1)} &:&R=2,\text{ }r=2t+1\geq 7,\text{ }q=5,\text{ }
r\neq 9,\text{ }\overline{\mu }_{5}(2)\approx \frac{8}{5}\text{,}  \notag \\
n \,&=\,&\sqrt{5}\cdot 5^{\frac{r-2}{2}}+\notag\\&&\left\{
\begin{array}{@{}l@{\,\,}l}
\overline{\ell }_{5}(\frac{r-1}{2},2)\smallskip & \text{if }r=4c+3 \\
(\overline{\ell }_{5}(\frac{r-1}{4},2)+\frac{1}{4})\cdot 5^{\frac{r-1}{4}}-
\frac{1}{4}\smallskip & \text{if }r=8c+5 \\
(\overline{\ell }_{5}(\frac{r-5}{4},2)+\frac{1}{2})\cdot 5^{\frac{r+3}{4}}-
\frac{1}{2} & \text{if }r=8c+1
\end{array}
\right. .  \label{form6_R=2_r=2t+1_q=5}
\end{eqnarray}

Now we construct infinite code families by using the $q^{m}
$-concatenating constructions in~\cite{DavNBCR2}. Terminology
and notation of \cite{DavNBCR2} will be used; in particular, we
are going to consider $2^{E}$-partitions, $ 2^{+}$-partitions,
and their cardinalities $h^{E}(H)$ and $h^{+}(H),$ see
\cite[Def.\thinspace 1,Rem.\thinspace 1]{DavNBCR2}. The
starting codes will be the codes associated to the 1-saturating
sets described in the part \emph{ C} of this section.

In \cite[Ex. 6, form. (33)]{DavNBCR2} an infinite family of $[n,n-r]_{q}2$
codes is constructed with
\begin{eqnarray}
\mathcal{A}_{2,q}^{(1)}\, &:&\,R=2,\text{ }r=2t+1\geq 3,\text{ }
q=(q^{\prime })^{2}\geq 16,\notag\\n\,&=\,&\left( 3-\frac{1}{\sqrt{q}}\right) q^{
\frac{r-2}{2}}+\lfloor q^{\frac{r-5}{2}}\rfloor ,  \notag \\
\overline{\mu }_{q}(2) &\,<\,&4.5-\frac{3}{\sqrt{q}}-\frac{17}{2q}+\frac{9}{q
\sqrt{q}}+\frac{5}{2q^{2}}.  \label{form6_R=2_r=2t+1_q=q'^2}
\end{eqnarray}
The starting code $V_{0}$ (denoted as $\mathcal{W}$) is based
on the previously mentioned 1-saturating $(3\sqrt{q}-1)$-set.
In \cite{DavNBCR2} it is noted that $h^{E}(H_{V_{0}})\leq 4$
and that this inequality allows us to obtain an effective
iterative code chain. A similar situation arises if one takes
as $V_{0}$ the $[n_{0}=2\sqrt{q}+2\sqrt[4]{q}+2,n_{0}-3]_{q}2$
code based on Theorem~\ref{Th4_q4-1sat}. We partition the
column set of the parity-check matrix into subsets
$T_{1},\ldots ,T_{4}$ so that $ |T_{1}|=|T_{3}|=2,$ $T_{1}\cup
T_{2}=\pi _{1},$ $T_{3}\cup T_{4}=\pi _{2},$ where $\pi
_{1},\pi _{2}$ are the disjoint Baer subplanes in
$PG(2,\sqrt{q})$ . An arbitrary point of $PG(2,q)\smallsetminus
\{\pi _{1}\cup \pi _{2}\}$ lies on a line through two points
belonging to the distinct subplanes. So, we obtain a
2-partition, see \cite[Def.\thinspace 1]{DavNBCR2} and
Definition~\ref{Def2_R,l code}. Moreover, as every point of a
subplane $\pi$ is a linear combination of two other points of
$\pi$, this $2$-partition is a 2$^{E}$-partition
\cite[Rem.\thinspace 1]{DavNBCR2} and $ h^{E}(H_{V_{0}})\leq
4.$ Now, by changing $3\sqrt{q}-1$ by $2\sqrt{q}+2\sqrt[
4]{q}+2$ in (\ref{form6_R=2_r=2t+1_q=q'^2}), we obtain the
following theorem.

\begin{theorem}
\label{Th6_R=2 q=(q')^4}For $q=(q^{\prime })^{4}$ there is an infinite
family of $[n,n-r]_{q}2$ codes with
\begin{eqnarray}
\mathcal{A}_{2,q}^{(1)}&:&R=2,\text{ }r=2t+1\geq 3,\text{ }
q=(q^{\prime })^{4},\notag\\n\,&=\,&\left( 2+\frac{2}{\sqrt[4]{q}}+\frac{2}{\sqrt{q
}}\right) q^{\frac{r-2}{2}}+\left\lfloor q^{\frac{r-5}{2}}\right\rfloor ,
\notag \\
\overline{\mu }_{q}(2) &\,<\,&2+\frac{4}{\sqrt[4]{q}}+\frac{6}{\sqrt{q}}+\frac{4
}{\sqrt[4]{q^{3}}}-\frac{2}{q}-\frac{8}{q\sqrt[4]{q}}.
\label{form6_R=2_q=(q')^4_inf-fam}
\end{eqnarray}
\end{theorem}

\begin{theorem}
\label{Th6_R=2_r=2t+1_nq<q}Let $q\geq 7.$ Assume that there
exists an $ [n_{q},n_{q}-3]_{q}2$ code $V_{0}$ with $n_{q}<q.$
Then there exists an infinite family of $[n,n-r]_{q}2$ codes
with
\begin{eqnarray}
\mathcal{A}_{2,q}^{(1)}&:&R=2,\,r=2t+1\geq 3,\,
q\geq 7,\,r\neq 9,13,\,a_{q}=\frac{n_{q}}{\sqrt{q}},  \notag \\
n &\,=\,&a_{q}\cdot q^{\frac{r-2}{2}}+2\lfloor q^{\frac{r-5}{2}}\rfloor +\left\{
\begin{array}{@{\,\,}c@{\,\,}l}
\renewcommand{\arraystretch}{1.2}
0 & \text{if }2p_{0}\leq q+1 \\
\lfloor q^{\frac{r-7}{2}}\rfloor & \text{if }2p_{0}>q+1
\end{array}
\right. ,  \notag \\
\overline{\mu }_{q}(2) &\approx &\frac{a_{q}^{2}}{2}-\frac{a_{q}^{2}}{q}+
\frac{2a_{q}}{q\sqrt{q}}.\label{form6_R=2_r=2t+1_nq<q}
\end{eqnarray}
For $r=9,13,$ $n=a_{q}\cdot
q^{\frac{r-2}{2}}+2q^{\frac{r-5}{2}}+q^{\frac{r-7
}{2}}+q^{\frac{r-9}{2}}$ holds.
\end{theorem}

\begin{IEEEproof}
Take $V_{0}$ as the starting code for the constructions of
\cite{DavNBCR2}. Then, changing $n_{q}$ by $p_{0},$ we use the
same argument of \cite[ Ex.\thinspace 6]{DavNBCR2} on partition
cardinalities $h^{+}(H_{\mathcal{W} }),$
$h^{E}(H_{\mathcal{W}})$. As a result,
(\ref{form6_R=2_r=2t+1_nq<q}) is obtained, cf.
\cite[form.\thinspace (32)]{DavNBCR2}.
\end{IEEEproof}

Theorem \ref{Th6_R=2_r=2t+1_nq<q} is the main tool to obtain infinite code
families with growing odd codimension.

\begin{theorem}
\label{Th6_R=2 q=(q')^6}For $q=(q^{\prime })^{6}$ there is an infinite
family of $[n,n-r]_{q}2$ codes with
\begin{eqnarray}
\mathcal{A}_{2,q}^{(1)}&:&R=2,\text{ }r=2t+1\geq 3,\text{ }
r\neq 9,13, \label{form6_R=2_r=2t+1_q=q'^6}\\
q&\,=\,&(q^{\prime })^{6},\text{ }q^{\prime }\text{ prime},\text{ }q^{\prime }\leq
73,\text{ }\notag \\
n &\,=\,&\left( 2+\frac{2}{\sqrt[6]{q}}+\frac{2}{\sqrt[3]{q}}+\frac{2}{\sqrt{q}}
\right) q^{\frac{r-2}{2}}+2\lfloor q^{\frac{r-5}{2}}\rfloor ,
  \notag \\
\overline{\mu }_{q}(2) &\,<\,&2+\frac{4}{\sqrt[6]{q}}+\frac{6}{\sqrt[3]{q}}+
\frac{8}{\sqrt{q}}+\frac{6}{\sqrt[3]{q^{2}}}+\frac{5}{\sqrt[6]{q^{5}}}.
\notag
\end{eqnarray}
\end{theorem}

\begin{IEEEproof}
The assertion follows from Theorem \ref{Th6_R=2_r=2t+1_nq<q}
and Corollary \ref{Cor4_q6-1sat}.
\end{IEEEproof}

\begin{lemma}
\label{Lem6_R=2_r=3}For an $[n_{q},n_{q}-3,3]_{q}2$ code
$V_{0}$ we have $ p_{0}\leq n_{q}-1$.
\end{lemma}

\begin{IEEEproof}
In a parity-check matrix $H$ of $V_{0}$ there are three linear
dependent columns. Let two of these columns form one subset of
a partition $\mathcal{P} _{0}$ of $H$, while the other subsets
of $\mathcal{P}_{0}$ contain precisely one column. By
Definition \ref{Def2_R,l code}, $\mathcal{P}_{0}$ is a
2-partition.
\end{IEEEproof}

\begin{theorem}
For any $q\leq 1217$, there exists an infinite family of $[n,n-r]_{q}2$
codes with
\begin{eqnarray}
\mathcal{A}_{2,q}^{(1)} &:&R=2,\text{ }r=2t+1\geq 3,\text{ }q\leq 1217,\label{form6_R=2_r=2t+1_byTable}
\\
&&r\neq 9,13,\text{ }a_{q}=\frac{\overline{\ell }_{q}(3,2)}{\sqrt{q}},\text{ }
\overline{\mu }_{q}(2)<\frac{a_{q}^{2}}{2},  \notag \\
n &\,=\,&a_{q}\cdot q^{\frac{r-2}{2}}+2\lfloor q^{\frac{r-5}{2}}\rfloor +\left\{
\begin{array}{c@{\,\,}l}
0 & \text{if }16\leq q\leq 1217 \\
\lfloor q^{\frac{r-7}{2}}\rfloor & \text{if }7\leq q\leq 13\smallskip
\end{array}
\right. ,   \notag \\
&&\overline{\mu }_{q}(2) <4.5\text{ if }q\leq 109,\text{ }\overline{\mu }
_{q}(2)<6.125\text{ if }q\leq 349,\notag \\&&\overline{\mu }_{q}(2)<8\text{ if }
q\leq 1217.  \notag
\end{eqnarray}
\end{theorem}

\begin{IEEEproof}
By Lemma \ref{Lem6_R=2_r=3} and Table \ref{table1}, for
$q=16,17$ we have $2p_{0}\leq q+1.$ Then the assertion follows
from Table \ref{table1} and Theorems~\ref {Th5_1sat_byTable}
and \ref{Th6_R=2_r=2t+1_nq<q}.
\end{IEEEproof}

For each of the infinite families
(\ref{form6_R=2_r=2t+1_q=3})-(\ref {form6_R=2_r=2t+1_byTable})
the covering density is bounded from above by a constant. If in
(\ref{form6_R=2_r=2t+1_nq<q}) we take as $V_{0}$ a code with
length $n_{q}\sim f(q)\sqrt{q}$, where $f(q)$ is some
increasing function of $q,$ such as in
(\ref{form6_R=2_r=2t+1_q^2/3}), then the asymptotic covering
density increases like $f^{2}(q).$ However for concrete $q$ new
code families can be supportable, see e.g. Table~\ref{table2}.

We end this section with Tables \ref{table5} and~\ref{table6},
which have been obtained from (\ref
{form6_DS}),(\ref{form6_l3(r+1,2)})-(\ref{form6_R=2_r=2t_3}),(\ref
{form6_R=2_r=2t+1_q=3})-(\ref{form6_R=2_r=2t+1_q=5}),(\ref
{form6_R=2_r=2t+1_byTable}), Table~\ref{table1},
Proposition~\ref{prop6_l4(4,2)}, and \cite[Tab.\thinspace
II]{BV}, \cite[Tab.\thinspace 1]{DavO1},\cite[ Tab.\thinspace
I]{DavO3},\cite[Tab.\thinspace I]{DavO5}.

\begin{table*}
 \caption{Upper Bounds $\overline{\ell }_{q}(r,2)$ on the
Length Function $\ell _{q}(r,2),$ $q=3,4,5,7,$ $r\leq 24$}
\renewcommand{\arraystretch}{1.0}
\centering
\begin{tabular}{@{}rrrrr|rrrrr@{}}
\hline
$r$ & $\overline{\ell }_{3}(r,2)$ & $\overline{\ell }_{4}(r,2)$ & $\overline{
\ell }_{5}(r,2)$ & $\overline{\ell }_{7}(r,2)$ & $r$ & $\overline{\ell }
_{3}(r,2)$ & $\overline{\ell }_{4}(r,2)$ & $\overline{\ell }_{5}(r,2)$ & $
\overline{\ell }_{7}(r,2) \phantom{\overline{\overline{H}}}$\\ \hline
$3$ & $4\centerdot $ & $5\centerdot $ & $6\centerdot $ & $6\centerdot $ & $
14 $ & $1822$ & $9522$ & $35000$ & $252105$ \\
$4$ & $8\centerdot $ & $9\centerdot $ & $11\phantom{\centerdot}$ & $15
\phantom{\centerdot}$ & $15$ & $2915$ & $19456$ & $78256$ & $741909$ \\
$5$ & $11\centerdot $ & $19\phantom{\centerdot}$ & $28\phantom{\centerdot}$
& $44\phantom{\centerdot}$ & $16$ & $5588$ & $37888$ & $175000$ & $1764735$
\\
$6$ & $22\phantom{\centerdot}$ & $37\phantom{\centerdot}$ & $56
\phantom{\centerdot}$ & $105\phantom{\centerdot}$ & $17$ & $8201$ & $77824$
& $410937$ & $5193363$ \\
$7$ & $40\phantom{\centerdot}$ & $85\phantom{\centerdot}$ & $131
\phantom{\centerdot}$ & $309\phantom{\centerdot}$ & $18$ & $16402$ & $151552$
& $875000$ & $12353145$ \\
$8$ & $76\phantom{\centerdot}$ & $154\phantom{\centerdot}$ & $281
\phantom{\centerdot}$ & $743\phantom{\centerdot}$ & $19$ & $24785$ & $316672$
& $1953828$ & $36353541$ \\
$9$ & $101\phantom{\centerdot}$ & $304\phantom{\centerdot}$ & $703
\phantom{\centerdot}$ & $2164\phantom{\centerdot}$ & $20$ & $49328$ & $
611328 $ & $4375000$ & $86472015$ \\
$10$ & $202\phantom{\centerdot}$ & $592\phantom{\centerdot}$ & $1400
\phantom{\centerdot}$ & $5145\phantom{\centerdot}$ & $21$ & $73811$ & $
1245184$ & $9853906$ & $254474787$ \\
$11$ & $323\phantom{\centerdot}$ & $1237\phantom{\centerdot}$ & $3153
\phantom{\centerdot}$ & $15141\phantom{\centerdot}$ & $22$ & $147622$ & $
2424832$ & $21875000$ & $605304105$ \\
$12$ & $620\phantom{\centerdot}$ & $2389\phantom{\centerdot}$ & $7031
\phantom{\centerdot}$ & $36407\phantom{\centerdot}$ & $23$ & $223073$ & $
4980736$ & $48831278$ & $1781323509$ \\
$13$ & $911\phantom{\centerdot}$ & $4948\phantom{\centerdot}$ & $16406
\phantom{\centerdot}$ & $106036\phantom{\centerdot}$ & $24$ & $443960$ & $
9699328$ & $109375000$ & $4237128735$ \\ \hline
\end{tabular}
\label{table5}
\end{table*}

\begin{table}
 \caption{Covering Densities $\overline{\mu
}_{q}(2,\mathcal{A} _{2,q}^{(\gamma )})$ of Infinite Families
$\mathcal{A}_{2,q}^{(\gamma )}$}
\renewcommand{\arraystretch}{1.0}
\centering
\begin{tabular}{@{}lcc|lcc|ccc@{}}
\hline
$q$ &$\gamma= 0$ &$ \gamma=1$ &$ q$ &$ \gamma=0$ & $\gamma=1$ &$ q$
&$\gamma= 0 $&$\gamma= 1 $\\ \hline
3 & 1.389 & 1.042 & 7 & 1.687 & 2.087 & 11 & 1.807 & 1.943 \\
4 & 1.504 & 1.587 & 8 & 1.729 & 1.880 & 13 & 1.838 & 2.183 \\
5 & 1.606 & 1.600 & 9 & 1.782 & 1.707 & 16 & 1.870 & 2.287 \\ \hline
\end{tabular}
\label{table6}
\end{table}

\section{\label{Sec7_R=3}Codes with Covering radius $R=3$}

\subsection{ Infinite Code Families of Codimension $r=3t$}

Let $q\geq 4.$ The geometrical construction (named
\textquotedblleft two ovals plus line\textquotedblright )
\cite[Th.\thinspace 7]{DavO4} gives a $ [3q+1,3q-5]_{q}3$ code.
So,
\begin{equation}
\ell _{q}(6,3)\leq 3q+1\text{ if }q\geq 4.  \label{form7_l_q(6,3)-bound}
\end{equation}

To our knowledge, no examples of $[n,n-6]_{q}3$ code with $n<3q+1$ are known.

\textbf{Open problem.} To prove that $\ell _{q}(6,3)=3q+1$ for
$q\geq 4.$\smallskip

The parity-check matrix of the code of \cite[Th. 7]{DavO4} is
modified in \cite[Th.\thinspace 6]{DavO5} and then it is used
as the starting point in $ q^{m}$-concatenating constructions.
As a result, an infinite family of $ [n,n-r]_{q}3$ codes is
obtained with the following parameters
\begin{eqnarray}
\mathcal{A}_{3,q}^{(0)} &:&R=3,\text{ }r=3t\geq 6,\text{ }q\geq 5,\text{ }
n=3q^{\frac{r-3}{3}}+q^{\frac{r-6}{3}},  \notag \\
r&\,\neq&\, 9,\text{ }\overline{\mu }_{q}(3) <\frac{9}{2}-\frac{9}{q}+\frac{3}{2q^{2}}+\frac{14}{
3q^{3}}-\frac{1}{2q^{4}}.  \label{form7_R=3_r=3t_DO5}
\end{eqnarray}
Also, in \cite{DavO5} it is shown that codes with parameters as
in~(\ref {form7_R=3_r=3t_DO5}) exist for $r=9$ if $q=16$ or
$q\geq 23.$ For $ q=7,8,11,13,17,19,$ DS of the codes
(\ref{form6_R=2_r=2t_2}) and the $ [\theta _{3,q},\theta
_{3,q}-3]_{q}1$ Hamming codes gives $ [n=3q^{2}+2q+1,n-9]_{q}3$
codes. For $q=5,9$ and $r=9,$ by (\ref {form6_R=2_r=2t_3}),
codes with length $n=3q^{2}+2q+2$ are obtained.

\subsection{ Infinite Code Families of Codimension
$r=3t+1$}

Let $q=3.$ DS of the codes of (\ref{form6_R=2_r=2t+1_q=3}) and
the $[\theta _{t,3},\theta _{t,3}-t]_{3}1$ Hamming codes forms
an infinite family of $ [n,n-r]_{3}3$ codes with
\begin{eqnarray}
\mathcal{A}_{3,3}^{(1)} &:&R=3,\text{ }r=3t+1\geq 7,\text{ }q=3,\text{ }
r\neq 10,\label{form7_R=3_r=3t+1_q=3}\\
&&\overline{\mu }_{3}(3)\approx 2.382,\notag \\
n &\,=\,&\frac{7}{4}\cdot 3^{\frac{2}{3}}\cdot 3^{\frac{r-3}{3}}+\left\{
\begin{array}{l@{\,\,}l}
-\frac{3}{4} & \text{if }r=6c+1 \\
\frac{3}{4}\cdot 3^{\frac{r+2}{6}}-\frac{3}{2}\smallskip & \text{if }r=12c+4
\\
\frac{3}{4}\cdot 3^{\frac{r+8}{6}}-\frac{3}{2} & \text{if }r=12c+10
\end{array}
\right. .  \notag
\end{eqnarray}
Also, $[431,415]_{3}3$ and $[3887,3865]_{3}3$ codes are given
in \cite[ Tab.\thinspace I, form.\thinspace (37)]{Dav95}, and a
$[14,7]_{3}3$ code is obtained in \cite[Tab.\thinspace II]{BV}.

\begin{theorem}
\label{Th7_R=3_qTabIII}Denote by $Q_{3}$ the set of values of $q$ for which
there is an $[\overline{\ell }_{q}(4,3),\overline{\ell }_{q}(4,3)-4,3]_{q}3$
code with minimum distance $d=3.$ Then, for $7\leq q\leq 563,$ there is an
infinite family of $[n,n-r]_{q}3$ codes with
\begin{eqnarray}
\mathcal{A}_{3,q}^{(1)} &:&R=3,\text{ }r=3t+1\geq 4,\text{ }7\leq q\leq 563,
\label{form7_R=3_TabIII}\\
&&b_{q}=\frac{\overline{\ell }_{q}(4,3)}{\sqrt[3]{q}},\text{ }
\overline{\mu }_{q}(3)<\frac{b_{q}^{3}}{6},\notag  \\
n &\,=&\,b_{q}\cdot q^{\frac{r-3}{3}}+\frac{q^{\frac{r-4}{3}}-1}{q-1}+\left\{
\begin{array}{c@{\,\,\,}c}
\frac{q^{\frac{r-4}{3}}-1}{q-1}\medskip & \text{if }q\in Q_{3} \\
\overline{\ell }_{q}(2\frac{r-4}{3},2) & \text{if }q\notin Q_{3}\medskip
\end{array}
\right. \text{,}   \notag \\
&&\overline{\mu }_{q}(3) <10.7\text{ if }q\leq 83,\text{ }\overline{\mu }
_{q}(2)<15.2\text{ if }q\leq 343,\notag\\&&\overline{\mu }_{q}(2)<20.9\text{
if }q\leq 563.  \notag
\end{eqnarray}
\end{theorem}

\begin{IEEEproof}
By Lemma \ref{Lem2_d-and-l}, for $q\in Q_{3}$ we have
$[\overline{\ell } _{q}(4,3),\overline{\ell
}_{q}(4,3)-4,3]_{q}3,\ell _{0}$ codes with $\ell _{0}\geq 1.$
By Table~\ref{table3}, for $q\geq 7$ we have that
$\overline{\ell } _{q}(4,3)\leq q$ if $q\in Q_{3}$ and
$\overline{\ell }_{q}(4,3)\leq q-2$ if $ q\notin Q_{3}.$ We
take the $[\overline{\ell }_{q}(4,3),\overline{\ell }
_{q}(4,3)-4]_{q}3$ codes of Table \ref{table3} as the codes
$V_{0}$ for Constructions QM$_{2}$ (if $q\in Q_{3}$) and
QM$_{4}$ (if $q\notin Q_{3}$), using the trivial partition and
letting $m\geq 1$. Now the assertion follows from (\ref
{form2_QM_2}) and (\ref{form2_QM_4}).
\end{IEEEproof}

We denote by $p^{(\ell )}(V)$ the upper bound of the minimal
possible cardinality of an $(R,\ell)$-partition for a
parity-check matrix of an $ [n,n-r]_{q}R,\ell $ code $V$.

\begin{theorem}
\label{Th7_R=3 q=g^3}For $q=(q^{\prime })^{3}\geq 64$ there exists an
infinite family of $[n,n-r]_{q}3$ codes with
\begin{eqnarray}
\mathcal{A}_{3,q}^{(1)} &:&R=3,\text{ }r=3t+1\geq 7,\text{ }
q=(q^{\prime })^{3}\geq 64,\label{form7_R=3_r=3t+1_q=q'^3}\\
n&\,=&\,\left( 4+\frac{4}{\sqrt[3]{q}}\right)
q^{\frac{r-3}{3}},\text{ }
\overline{\mu }_{q}(3) <\frac{32}{3}+\frac{32}{\sqrt[3]{q}}+\frac{32}{\sqrt
[3]{q^{2}}}-\frac{64}{3q}.  \notag
\end{eqnarray}
\end{theorem}

\begin{IEEEproof}
The 2-saturating $(4q^{\prime }+4)$-set $B$ of Corollary \ref
{Cor4_3fold_PG(3,q)} consists of pairwise skew lines of
$PG(3,q^{\prime })$. As $q^{\prime }\geq 4$, it can be shown
that the related code $C_{B}$ is a (3,3)-object,
see~(\ref{form4_R=3_r=3t+1_2sat-4lines}), Definition \ref
{Def2_R,l code}, Lemma \ref{Lem2_d-and-l}, and Remark \ref
{Rem4_points-on-line}. We take $C_{B}$ as the starting
$[n_{0}=4\sqrt[3]{q} +4,n_{0}-4,3]_{q}3,3$ code $V_{0}$ for
Construction QM$_{3}$ of Section~\ref {Sec2_q^m_concat}. The
trivial partition gives $p_{0}=p^{(3)}(V_{0})=n_{0}<q$ . So, we
take $m\geq 1$ and obtain a family of
$[n=q^{m}n_{0},n-(4+3m)]_{q}3$ codes.
\end{IEEEproof}

\subsection{ Infinite Code Families of Codimension
$r=3t+2$}

Let $q=3.$ DS of the codes of (\ref{form6_R=2_r=2t+1_q=3}) and the $[\theta
_{t+1,3},\theta _{t+1,3}-(t+1)]_{3}1$ Hamming codes forms an infinite family
of $[n,n-r]_{3}3$ codes with
\begin{eqnarray}
\mathcal{A}_{3,3}^{(2)} &:&R=3,\text{ }r=3t+2\geq 8,\text{ }q=3,\text{ }
r\neq 11, \label{form7_R=3_r=3t+2_q=3_DS}\\
&&\overline{\mu }_{3}(3)\approx 3.082, \notag\\
 n&\, =&\,\frac{11}{4}\sqrt[3]{3}\cdot 3^{\frac{r-3}{3}}+\left\{
\begin{array}{l@{\,\,\,}l}
-\frac{3}{4} & \text{if }r=6c+2 \\
\frac{3}{4}\cdot 3^{\frac{r+1}{6}}-\frac{3}{2}\smallskip & \text{if }r=12c+5
\\
\frac{3}{4}\cdot 3^{\frac{r+7}{6}}-\frac{3}{2} & \text{if }r=12c+11
\end{array}
\right. .  \notag
\end{eqnarray}

Also, $[674,657]_{3}3$ and $[6074,6051]_{3}3$ codes are given
in \cite[ Tab.\thinspace I, form.\thinspace (38)]{Dav95}.

\begin{theorem}
\label{Th7_R=3_qTabIV} For $3\leq q\leq 43,$ there exists an infinite family
of $[n,n-r]_{q}3$ codes with
\begin{eqnarray}
\mathcal{A}_{3,q}^{(2)} &:&R=3,\text{ }r=3t+2\geq 5,\text{ }r\neq 8,\text{ }
3\leq q\leq 43, \label{form7_R=3_TabIV}\\
&&c_{q}=\frac{\overline{\ell }_{q}(5,3)}{\sqrt[3]{q^{2}}
},\text{ }\overline{\mu }_{q}(3)<\frac{c_{q}^{3}}{6},  \notag \\
n &=&c_{q}\cdot q^{\frac{r-3}{3}}+\frac{q^{\frac{r-5}{3}}-1}{q-1}+\left\{
\begin{array}{c@{\,\,}c}
\frac{q^{\frac{r-5}{3}}-1}{q-1}\medskip & \text{if }q\neq 2,5,19 \\
\overline{\ell }_{q}(2\frac{r-5}{3},2) & \text{if }q=2,5,19\medskip
\end{array}
\right. \text{,}  \notag \\
&&\overline{\mu }_{q}(3) <10.7\text{ if }q\leq 27,\text{ }\overline{\mu }
_{q}(2)<12.4\text{ if }q\leq 32,\notag
\\&&\overline{\mu }_{q}(2)<20.9\text{ if
}q\leq 43.  \notag
\end{eqnarray}
\end{theorem}

\begin{IEEEproof}
By Lemma \ref{Lem2_d-and-l} and Table \ref{table4}, for $q\neq
2,5,19$ we have $[ \overline{\ell }_{q}(5,3),\overline{\ell
}_{q}(5,3)-5,3]_{q}3,\ell _{0}$ codes with $\ell _{0}\geq 1.$
By Table \ref{table4}, for $q\geq 3$ we have that $
\overline{\ell }_{q}(5,3)\leq q^{2}.$ We take the
$[\overline{\ell } _{q}(5,3),\overline{\ell }_{q}(5,3)-5]_{q}3$
codes of Table \ref{table4} as the codes $ V_{0}$ for
Constructions QM$_{2}$ (if $q\neq 2,5,19)$ and QM$_{4}$
(otherwise) using the trivial partition and letting $m\geq 2$.
Now the assertion follows from (\ref{form2_QM_2}) and
(\ref{form2_QM_4}).
\end{IEEEproof}

\begin{theorem}
\label{Th7_R=3 q=g^3 r=2t+2} For $q=(q^{\prime })^{3}\geq 27$ there exists
an infinite family of $[n,n-r]_{q}3$ codes with
\begin{eqnarray}\mathcal{A}_{3,q}^{(2)} :R=3,\text{ }r=3t+2\geq 8,\text{ }q=(q^{\prime
})^{3}\geq 27,\qquad\label{form7_fam_R=3_r=3t+2_q=q'^3}\\
n=\left( 9-\frac{8}{\sqrt[3]{q}}+\frac{4}{\sqrt[3]{
q^{2}}}\right) q^{\frac{r-3}{3}},\,\overline{\mu }_{q}(3) <\frac{243}{2}-\frac{324}{\sqrt[3]{q}}+\frac{72}{
\sqrt[3]{q^{2}}}. \notag
\end{eqnarray}
\end{theorem}

\begin{IEEEproof}
Let $S$ be as (\ref{form4_R=3_r=3t+2_2sat-9planes}). For any
plane of $S$, let
$\{P_{1},P_{2}\},\{P_{3},P_{4}\},\{P_{5},\ldots ,P_{(q^{\prime
})^{2}+q^{\prime }+1}\}$ be a partition of the set of its
points such that $ P_{1},P_{2}\notin l_{3,4}$ and
$P_{3},P_{4}\notin l_{1,2}$, where $l_{i,j}$ is the line
through points $P_{i},P_{j}$. It can be easily shown that if $
u\in \{2,3\}$, then every point of the plane is equal to a
linear combination with nonzero coefficients of $u$ other
points belonging to distinct subsets of the partition. The
corresponding partition of the columns of the parity-check
matrix of the related code $\mathcal{C}_{S}$ is a
(3,3)-partition with $p^{(3)}(\mathcal{C} _{S})=3\cdot 9=27\leq
q.$ Therefore we may take $\mathcal{C}_{S}$ as the starting
$[n_{0}=9\sqrt[3]{q^{2}}-8\sqrt[3]{q}+4,n_{0}-4]_{q}3,3$ code $
V_{0} $ for Construction QM$_{3}$ with $m\geq 1$.
\end{IEEEproof}

\section{\label{Sec8_R>=4}Codes with Covering Radius $R\geq 4$}

\subsection{ Infinite Code Families of Codimension $r=Rt$
and Arbitrary $ q$}

In this Section we obtain a code $V$ of covering radius $R\geq 4$ and
codimension $Rt$ from DS of $g_{2}$ codes $V_{2}$ with radius two and $g_{3}$
codes $V_{3}$ with radius three. More precisely, let
\begin{equation}
V=\underbrace{V_{2}\oplus \ldots \oplus V_{2}}_{g_{2}\text{ times}}\oplus
\underbrace{V_{3}\oplus \ldots \oplus V_{3}}_{g_{3}\text{ times}}
\label{form8_DS_r=Rt}
\end{equation}
where $V$ is an $[n,n-Rt]_{q}R$ code, $V_{2}$ is an $[n_{2},n_{2}-2t]_{q}2$
code, $V_{3}$ is an $[n_{2},n_{2}-3t]_{q}3$ code, $n=g_{2}n_{2}+g_{3}n_{3},$
$2g_{2}+3g_{3}=R,$ and
\begin{equation}
g_{2}=\left\{
\begin{array}{l@{\,\,\,}l}
0 & \text{if }R\equiv 0\pmod 3 \\
1 & \text{if }R\equiv 2\pmod 3 \\
2 & \text{if }R\equiv 1\pmod 3
\end{array}
\right. ,\text{ }g_{3}=\left\lceil \frac{R}{3}\right\rceil -g_{2}.
\label{form8_(g2,g3)}
\end{equation}

\begin{theorem}
\label{Th8_t=2}Let $R\geq 4$ and let $q\geq 4.$ Then there
exists an $ [n=Rq+\left\lceil R/3\right\rceil
,n-2R,3]_{q}R,\ell $ code with $\ell \geq 1.$
\end{theorem}

\begin{IEEEproof}
Geometrical constructions of a $[2q+1,2q-3]_{q}2$ code $V_{2}$
(\textquotedblleft oval plus line\textquotedblright ) and of a
$ [3q+1,3q-5]_{q}3 $ code $V_{3}$ (\textquotedblleft two ovals
plus line\textquotedblright ) are given in\ \cite[p.\thinspace
104]{BrPlWi}$,$ \cite[Th.\thinspace
5.1]{Dav95},\cite[Th.7]{DavO4}. Using these codes in
(\ref{form8_DS_r=Rt}) and (\ref{form8_(g2,g3)}) with $t=2,$ we
obtain an $ [n=Rq+\left\lceil R/3\right\rceil ,n-2R]_{q}R$ code
$V$. Minimum distance $ d=3$ follows from the fact that the
point sets associated to $V_2$ and $V_3$ contain triples of
collinear points. The value $\ell \geq 1$ follows from Lemma
\ref{Lem2_d-and-l}.
\end{IEEEproof}

\textbf{Open problem.} To obtain $[n,n-2R]_{q}R$ codes with
$R\geq 4,$ $q\geq 4,$ $n<Rq+\left\lceil R/3\right\rceil .$ In
particular, for $R\geq 4,$ to generalize the geometrical
constructions \textquotedblleft oval plus
line\textquotedblright\ and \textquotedblleft two ovals plus
line\textquotedblright.

\begin{theorem}
There exist infinite families of $[n,n-r]_{q}R$ codes with the
parameters
\begin{eqnarray}
\text{i) }\mathcal{A}_{R,q}^{(0)}&:&R\geq 4,\text{ }r=Rt\geq
5R,\text{ }
q\geq 7,\text{ }q\neq 9,\notag\\
n&\,=&\,Rq^{\frac{r-R}{R}}+\left\lceil R/3
\right\rceil q^{\frac{r-2R}{R}},\text{ }r\neq 6R.
\label{form8_R>=4_r=tR_q-general}\\
\text{ii) }
\mathcal{A}_{R,q}^{(0)} &:&R\geq 4,\,r=Rt\geq 2R,
\,q=5,9,\,r\neq 3R,4R,6R, \notag \\
n &\,=\,&R q^{\frac{r-R}{R}}+\left( \left\lceil R/3\right\rceil +
g_{2}\cdot q^{-1}\right) q^{\frac{r-2R}{R}}.  \label{form8_R>=4_r=Rt_q=5,9}
\end{eqnarray}
\end{theorem}

\begin{IEEEproof}
We use the construction of (\ref{form8_DS_r=Rt}),
(\ref{form8_(g2,g3)}) with the codes $V_{2}$ and $V_{3}$ taken
from (\ref{form6_R=2_r=2t_2}),(\ref {form7_R=3_r=3t_DO5}) and
(\ref{form6_R=2_r=2t_3}),(\ref {form7_R=3_r=3t_DO5}) for the
cases i) and ii), respectively.
\end{IEEEproof}

It should be noted that the main term of the asymptotic
covering density $ \overline{\mu
}_{q}(R,\mathcal{A}_{R,q}^{(0)})$ for the family of (\ref
{form8_R>=4_r=tR_q-general}) is $\frac{R^{R}}{R!};$ it does not
depend on $ q. $

By the results on cases $r=8,12$ and $r=9$ reported after (\ref
{form6_R=2_r=2t_2}),(\ref{form6_R=2_r=2t_3}), and
(\ref{form7_R=3_r=3t_DO5}), one can easy fill up gaps in
(\ref{form8_R>=4_r=tR_q-general}),(\ref{form8_R>=4_r=Rt_q=5,9})
for codes with $r=3R,4R,$ and $6R$.

\subsection{Infinite Code Families of Codimension $r=Rt+1,$ \\$
q=(q^{\prime })^{R}$}

\begin{theorem}
\label{Th8_r=R+1q=g^R+1}Let $q=(q^{\prime })^{R}.$ Then there exists an
infinite family of $[n,n-r]_{q}R$ codes with
\begin{eqnarray}
\mathcal{A}_{R,q}^{(1)}\text{ } &:&\text{ }R\geq 4,\text{ }r=Rt+1,\text{ }
q=(q^{\prime })^{R},\label{form8_R>=4_r=Rt+1}\\
&&t=1\text{ and }t\geq t_{0},\text{ }q^{t_{0}-1}\geq n_{R,q}^{(1)},  \notag \\
n_{R,q}^{(1)} &=&(\sqrt[R]{q}-1)\left( \frac{R(R+1)}{2}-2\right) +R+5,
 \notag\\
n &\,=\,&n_{R,q}^{(1)}\cdot q^{\frac{r-(R+1)}{R}}+\notag\\
&&\left\{
\begin{array}{l@{\,\,}l}
0 & \text{if }q^{\prime }\geq 4 \\
w\frac{q^{\frac{r-(R+1)}{R}}-1}{q-1},\text{ }w\in \{0,1\}, & \text{if }
q^{\prime }=3
\end{array}
\right. .  \notag
\end{eqnarray}
\end{theorem}

\begin{IEEEproof}
As the starting code $V_{0}$ for Constructions QM$_{2}$,
QM$_{3}$ we take an
$[n_{R,q}^{(1)},n_{R,q}^{(1)}-(R+1),3]_{q}R$ code
$\mathcal{C}_{K}$ related to the $(R-1)$-saturating set
$K\subset PG(R,q^{\prime })\subset PG(R,q)$ described in
Corollary~\ref{Cor4_PG(ro+1,q)_ro-sat}, see also (\ref
{form4_inductive_strong}), Construction A and Corollary
\ref{Cor4_Nfold_indu} . Note that $K$ contains four pairwise
skew lines of $PG(R,q^{\prime })$, whereas for other
$\frac{R(R+1)}{2}-6\geq 2R-4$ all but one point belong to $ K$.
These latter lines are partitioned into $R-3$ sets of
concurrent lines. By Definition~\ref{Def2_R,l code} and Remark
\ref{Rem4_points-on-line}, the code $\mathcal{C}_{K}$ is an
$(R,\ell _{0})$-object with $\ell _{0}=R$ if $ q\geq 4$ and
$\ell _{0}\geq R-1$ if $q=3.$ The trivial partition of its
parity-check matrix is an $(R,\ell _{0})$-partition into
$n_{R,q}^{(1)}\leq q^{t_{0}-1}$ subsets. Finally, we use
(\ref{form2_QM_2}) and (\ref {form2_QM_3}) to get the
assertion.
\end{IEEEproof}

It should be noted that the main term of the asymptotic
covering density $ \overline{\mu
}_{q}(R,\mathcal{A}_{R,q}^{(1)})$ for the family of (\ref
{form8_R>=4_r=Rt+1})\ is $\frac{(R^{2}+R)^{R}}{2^{R}R!};$ it
does not depend on $q.$

\subsection{Infinite Code Families of Codimension\\
$r=Rt+2,\ldots ,R(t+1)-1,$ $q=(q^{\prime })^{R}$}

To our knowledge, for $R\geq 4,$ $r=Rt+2,\ldots ,R(t+1)-1,$ no infinite
families with density asymptotically independent on $q$ are known.

\begin{theorem}
\label{Th8_r=R+a}Let $q=(q^{\prime })^{R}.$ Fix $\gamma \in \{2,3,\ldots
,R-1\}.$ Then there exists an infinite family of $[n,n-r]_{q}R$ codes with
\begin{eqnarray}
\mathcal{A}_{R,q}^{(\gamma )}&:&R\geq 4,\text{ }r=Rt+\gamma
,\text{ }q=(q^{\prime })^{R}, \label{form8_R>=4_r=Rt+h} \\
\gamma &\,=&\,2,3,\ldots ,R-1,\text{ }
t =1\text{ and }t\geq t_{0},\text{ }q^{t_{0}-1}\geq n_{R,q}^{(\gamma )},
\text{ }  \notag \\
n_{R,q}^{(\gamma )} &\,=\,&\frac{\sum_{i=1}^{\gamma +1}(\sqrt[R]{q}-1)^{i} {
\binom{{R+\gamma} }{i}}}{\sqrt[R]{q}-1}\sim \binom{R+\gamma }{R-1}q^{\frac{
\gamma }{R}},   \notag \\
n &\,=\,&n_{R,q}^{(\gamma )}\cdot q^{\frac{r-(R+\gamma )}{R}}+w\frac{q^{\frac{
r-(R+\gamma )}{R}}-1}{q-1},\text{ }0\leq w\leq R-3.  \notag
\end{eqnarray}
\end{theorem}

\begin{IEEEproof}
As the starting code $V_{0}$ for Constructions QM$_{2}$,
QM$_{3}$ we take an $[\overline{n}_{R,q}^{(\gamma
)},\overline{n}_{R,q}^{(\gamma )}-(R+\gamma ),3]_{q}R$ code
$\mathcal{C}_{B_{R-1}}$ related to the $(R-1)$-saturating set
$B_{R-1}\subset PG(R+\gamma +1,q^{\prime })\subset $
$PG(R+\gamma +1,q)$ of Lemma \ref{Lem4_key} and
Theorem~\ref{Th4_Napoleon}. In (\ref{form4_B_k}
),(\ref{form4_ro-sat_general}) we put $k=\rho ,$ $v-\rho
=\gamma \geq 2,$ $ \rho =R-1.$ By Definition~\ref{Def2_R,l
code} and Remark \ref {Rem4_points-on-line}, the code
$\mathcal{C}_{B_{R-1}}$ is an $(R,\ell _{0})$-object with $\ell
_{0}\geq 3$ as the set $B_{R-1}$ contains lines. The trivial
partition of a parity-check matrix of $\mathcal{C}_{B_{R-1}}$
is an $ (R,\ell _{0})$-partition into $n_{R,q}\leq q^{t_{0}-1}$
subsets. Finally, we use (\ref{form2_QM_2}) and
(\ref{form2_QM_3}).
\end{IEEEproof}

It should be noted that the main term of the asymptotic
covering density $ \overline{\mu
}_{q}(R,\mathcal{A}_{R,q}^{(\gamma )})$ for the family of (\ref
{form8_R>=4_r=Rt+h})\ is $\left( \frac{(R+\gamma
)^{R-1}}{(R-1)!}\right) ^{R}\cdot \frac{1}{R!},$ which does not
depend on $q.$

\section{\label{Sec9_R nonprime}Codes with Nonprime Covering Radius $
R=sR^{\prime }$}

We consider the case when covering radius $R$ is nonprime, i.e.
$ R=sR^{\prime }$ with integer $s$ and $R^{\prime }.$

\begin{lemma}
\label{Lem9_dirsum}Let $R=sR^{\prime }.$ Assume that there
exists an $[n^{\prime },n^{\prime }-(R^{\prime }t+t^{\prime
})]_{q}R^{\prime }$ code $\mathcal{C}_{0}$ with $R^{\prime
}>t^{\prime }$. Then there exists an $[n^{\prime
}\frac{R}{R^{\prime }} ,n^{\prime }\frac{R}{R^{\prime
}}-(Rt+\frac{R}{R^{\prime }}t^{\prime })]_{q}R $ code
$\mathcal{C}$. Moreover, if the starting code $\mathcal{C}
_{0}$ is short the new code $\mathcal{C}$ is short too.
\end{lemma}

\begin{IEEEproof}
We apply Construction DS to $s$ copies of $\mathcal{C}_{0}$. If
the code $ \mathcal{C}_{0}$ is short then $n^{\prime
}=O(q^{(R^{\prime }t+t^{\prime }-R^{\prime })/R^{\prime }})$
or, in other words, $n^{\prime }=cq^{(R^{\prime }t+t^{\prime
}-R^{\prime })/R^{\prime }}$ where $c$ is a constant
independent of $q.$ Also, $(R^{\prime }t+t^{\prime }-R^{\prime
})/R^{\prime }=(Rt+st^{\prime }-R)/R.$ Therefore $n^{\prime
}\frac{R}{ R^{\prime }}=c\frac{R}{R^{\prime
}}q^{(Rt+\frac{R}{R^{\prime }}t^{\prime }-R)/R}.$ Then the
assertion is proved.
\end{IEEEproof}

\begin{corollary}
\label{cor9_evenR}For even $R\geq 4$ there exist infinite
families $\mathcal{ A}_{R,q}^{(R/2)}$ of $[n,n-r]_{q}R$ codes
with codimension $r=Rt+\frac{R}{2}$ and the following
parameters:
\begin{eqnarray}
\text{i) }q&\,=&\,(q^{\prime })^{2},\text{ }t\geq 1,\label{form9_evenR_q=q'^2}\\n&\,=&\,\frac{R}{2}\left( 3-
\frac{1}{\sqrt{q}}\right) q^{\frac{r-R}{R}}+\frac{R}{2}\left\lfloor \frac{1}{
\sqrt{q}}q^{\frac{r-2R}{R}}\right\rfloor .\notag\\
\text{ii) }q\,&=&\,(q^{\prime })^{4},\text{ }t\geq 1,\label{form9_evenR_q=q'^4}\\n&\,=&\,R\left( 1+\frac{1}{
\sqrt[4]{q}}+\frac{1}{\sqrt{q}}\right) q^{\frac{r-R}{R}}+\frac{R}{2}
\left\lfloor \frac{1}{\sqrt{q}}q^{\frac{r-2R}{R}}\right\rfloor .\notag\\
\text{iii) }q &\,=&\,(q^{\prime })^{6},\text{ }q^{\prime }\text{ prime},\text{ }
q^{\prime }\leq 73,\text{ }t\geq 1,\text{ }t\neq 4,6, \label{form9_evenR_q=q'^6}  \\
n &\,=&\,R\left( 1+\frac{1}{\sqrt[6]{q}}+\frac{1}{\sqrt[3]{q}}+\frac{1}{\sqrt{q}}
\right) q^{\frac{r-R}{R}}+R\left\lfloor \frac{1}{\sqrt{q}}q^{\frac{r-2R}{R}
}\right\rfloor .\notag
\end{eqnarray}
\end{corollary}
\begin{IEEEproof}
Put $R^{\prime }=2$ and use the codes of
(\ref{form6_R=2_r=2t+1_q=q'^2}),(
\ref{form6_R=2_q=(q')^4_inf-fam}), and
(\ref{form6_R=2_r=2t+1_q=q'^6}) as the code $\mathcal{C}_{0}$
of Lemma~\ref{Lem9_dirsum}.
\end{IEEEproof}

\begin{corollary}
\label{cor9_R=3u}Let $q=(q^{\prime })^{3}$ and assume that $3$ divides $R.$
Then there exist infinite families of $[n,n-r]_{q}R$ codes with
\begin{eqnarray}
\text{i) }\mathcal{A}_{R,q}^{(\frac{R}{3})}&:&R=3s,\text{ }
r=Rt+\frac{R}{3},\text{ }q=(q^{\prime })^{3}\geq 64,\notag\\t\,&\,\geq& 1,\text{ }
n=\frac{4R}{3}\left( 1+\frac{1}{\sqrt[3]{q}}\right) q^{\frac{r-R}{R}}.\label{form9_3|R_q=q'^3}\\
\text{ii) }\mathcal{A}_{R,q}^{(\frac{2R}{3})}&:&R=3s,\text{ }r=Rt+\frac{2R}{3},
\text{ }q=(q^{\prime })^{3}\geq 27,\notag\\t&\,\geq\,& 1,\text{ }n=\frac{R}{3}
\left( 9-\frac{8}{\sqrt[3]{q}}+\frac{4}{\sqrt[3]{q^{2}}}\right) q^{\frac{r-R
}{R}}.\label{form9_3|R_q=q'^3_2}
\end{eqnarray}
\end{corollary}

\begin{IEEEproof}
Put $R^{\prime }=3$ and use the codes of
(\ref{form7_R=3_r=3t+1_q=q'^3}) and
(\ref{form7_fam_R=3_r=3t+2_q=q'^3}) as the code
$\mathcal{C}_{0}$ of Lemma \ref{Lem9_dirsum}.
\end{IEEEproof}

\begin{corollary}
\label{cor9_R=s'R'}Let $R=sR^{\prime }.$ Let $q=(q^{\prime })^{R^{\prime }}.$
Then there exist an infinite family of $[n,n-r]_{q}R$ codes with
\begin{eqnarray}
\mathcal{A}_{R,q}^{(s)} &:&R=sR^{\prime },\text{ }R^{\prime }\geq 4,\text{ }
r=Rt+\frac{R}{R^{\prime }},\text{ }q=(q^{\prime })^{R^{\prime }},\notag \\
&&t=1
\text{ and }t\geq t_{0},\text{ }q^{t_{0}-1}\geq n_{R^{\prime },q}^{(1)},
\notag \\
n_{R^{\prime },q}^{(1)} &=&(\sqrt[R^{\prime }]{q}-1)\left( \frac{R^{\prime
}(R^{\prime }+1)}{2}-2\right) +R^{\prime }+5, \label{form9_R=s'R'} \\
n &=&\frac{R}{R^{\prime }}\cdot n_{R^{\prime },q}^{(1)}\cdot q^{\frac{r-(R+s)
}{R}}+\notag
\\&&\left\{
\begin{array}{l@{\,\,}l}
0 & \text{if }q^{\prime }\geq 4 \\
w\frac{R}{R^{\prime }}\cdot \frac{q^{\frac{r-(R+s)}{R}}-1}{q-1},\text{ }w\in
\{0,1\}, & \text{if }q^{\prime }=3
\end{array}
\right. .  \notag
\end{eqnarray}
\end{corollary}

\begin{IEEEproof}
We use the codes of (\ref{form8_R>=4_r=Rt+1}) as the code
$\mathcal{C}_{0}$ of Lemma~\ref{Lem9_dirsum}.
\end{IEEEproof}

\begin{corollary}
\label{cor9_R=s'R'_2}Let $R=sR^{\prime }.$ Let $q=(q^{\prime })^{R^{\prime
}}.$ Fix $\gamma \in \{2,3,\ldots ,R-1\}.$ Then there exists an infinite
family of $[n,n-r]_{q}R$ codes with
\begin{eqnarray}
\mathcal{A}_{R,q}^{(s)} &:&R=sR^{\prime },\text{ }R^{\prime }\geq 4,\text{ }
r=Rt+\frac{R}{R^{\prime }}\gamma ,\text{ }q=(q^{\prime })^{R^{\prime }},
 \notag \\
\gamma &\,=\,&2,3,\ldots ,R-1,\text{ } t= 1\text{ and }t\geq t_{0},\text{ }q^{t_{0}-1}\geq n_{R^{\prime
},q}^{(\gamma )},  \notag \\
n_{R^{\prime },q}^{(\gamma )} &=&\frac{\sum_{i=1}^{\gamma +1}(\sqrt[
R^{\prime }]{q}-1)^{i}\binom{R^{\prime }+\gamma }{i}}{\sqrt[R^{\prime }]{q}-1
}\sim \binom{R^{\prime }+\gamma }{R^{\prime }-1}q^{\frac{\gamma }{R^{\prime }
}},\text{ }  \label{form9_R=s'R'_general} \\
n &\,=\,&\frac{R}{R^{\prime }}\cdot n_{R^{\prime },q}^{(\gamma )}\cdot q^{\frac{
r-(R+s\gamma )}{R}}+w\frac{q^{\frac{r-(R+s\gamma )}{R}}-1}{q-1},\notag\\
&&0\leq w\leq R-3.  \notag
\end{eqnarray}
\end{corollary}

\begin{IEEEproof}
We use the codes of (\ref{form8_R>=4_r=Rt+h}) as the code
$\mathcal{C}_{0}$ of Lemma~\ref{Lem9_dirsum}.
\end{IEEEproof}

It should be noted that for the infinite families
(\ref{form9_evenR_q=q'^2})-(\ref{form9_R=s'R'_general}), the
main term of the lower limit of covering density $\overline{\mu
}_{q}(R,\mathcal{A} _{R,q}^{(\gamma )})$ is, respectively,
$\frac{R^{R}}{R!}\left( \frac{3}{2} \right) ^{R},$
$\frac{R^{R}}{R!},$ $\frac{R^{R}}{R!},$ $\frac{R^{R}}{R!}
\left( \frac{4}{3}\right) ^{R},\smallskip$
$\frac{R^{R}}{R!}3^{R},$ $\frac{R^{R}}{R!} \left(
\frac{R^{\prime }+1}{2}\right) ^{R},$ $\frac{R^{R}}{R!}\left(
\frac{(R^{\prime }+\gamma )^{R^{\prime }-1}}{R^{\prime
}!}\right)^{R}$. All these terms do not depend on $q.$

\begin{remark}
It should be emphasized that codes of Corollaries
\ref{cor9_evenR}-\ref {cor9_R=s'R'_2} are \textquotedblleft
short\textquotedblright\ for $ R=sR^{\prime }$ though as a rule
in these codes $q\neq (q^{\prime })^{R}.$ Usually we have this
property when $q=(q^{\prime })^{R}.$
\end{remark}

\section{\label{Sec10_Conclusion}Conclusion}

We considered infinite sequences $\mathcal{A}_{R,q}$ of linear
nonbinary covering codes ${\mathcal{C}}_{n}$ of type
$[n,n-r_{n}]_{q}R.$ Without loss of generality, we assumed that
the sequence of codimension $r_{n}$ is not decreasing. For a
given family $\mathcal{A}_{R,q},$ the covering radius $R$ and
the size $q$ of the underlying Galois field are fixed. We
considered also \emph{infinite sets of the families\
}$\mathcal{A}_{R,q}$, where $R$ is fixed but $q$ ranges over an
infinite set of prime powers.

Each infinite family $\mathcal{A}_{R,q}$ consists of
\emph{supporting }and \emph{filling }codes$.$ The supporting
codes are the codes ${\mathcal{C}}_n$ such that $r_n>r_{n+1}$.
Non-supporting codes are called filling codes. This terminology
is motivated by the fact that the parameters of the codes in a
family are completely determined by those of its supporting
codes. However, considering filling codes is necessary to
investigate not only the lower limit ($\liminf $) of the
covering densities of a family, but also its upper limit
($\limsup $).

Such lower and upper limits (denoted by $\overline{\mu
}_{q}(R,\mathcal{A} _{R,q})$ and $\mu _{q}^{\ast
}(R,\mathcal{A}_{R,q})$ respectively) are the most considerable
asymptotic features of families $\mathcal{A}_{R,q}$. It is also
relevant how these limits depend on $q$ in infinite sets of
families \emph{\ }$\mathcal{A}_{R,q}$ with fixed $R.$ We showed
that for the \emph{ upper limit} the best possibility is $\mu
_{q}^{\ast }(R,\mathcal{A} _{R,q})=O(q).$ The problem of
constructing infinite sets of families $ \mathcal{A}_{R,q}$
with $\mu _{q}^{\ast }(R,\mathcal{A}_{R,q})=O(q)$ is open in
the general case. We call it \emph{Open Problem~1}. In the
literature, a solution to Open Problem 1 was known only for
$R=2,$ $q$ square.

We first showed in Introduction that Open Problem 1 for
covering radius $R$ is solved provided that a solution to the
following \emph{\ Open Problem 2 } is achieved: construct $R$\
infinite code families $\mathcal{A}_{R,q}^{(\gamma)}$, $
\gamma=0,\ldots,R-1$, such that $\overline{\mu
}_{q}(R,\mathcal{A}_{R,q}^{(\gamma )})=O(1)$\ holds. Here
$\mathcal{A}_{R,q}^{(\gamma )}$ is an infinite family such that
its supporting codes are a sequence of
$[n_{u},n_{u}-r_{u}]_{q}R$ codes with codimension
$r_{u}=Ru+\gamma $ and length $n_{u}=f_{q}^{(\gamma )}(r_{u}),$
where $u\ge u_0$; $f_{q}^{(\gamma )}$ is an increasing function
for a fixed $q$.

The main achievement of the paper is a solution to Open Problem
2 (and, thereby, to Open Problem 1) for an arbitrary covering
radius $R\geq 2$. This solution consists of infinite sets of
families $\mathcal{A}_{R,q}$ where $ q=(q^{\prime })^{R}$,
$q^{\prime }$ is power of prime$.$ The main tool was using
codes related to saturating sets in projective spaces as
starting points for $q^{m}$-concatenating constructions of
covering codes. Combining $ q^{m}$-concatenating constructions
and the saturating sets turned out to be very effective.

In addition, the methods used for solving Open Problems 1 and 2
allowed us to obtain a number of results on covering codes of
independent interest. In particular, we obtained many new upper
bounds on the asymptotic covering density $\overline{\mu
}_{q}(R,\mathcal{A}_{R,q}^{(\gamma )})$ for distinct $R$ and $
\gamma .$ We obtained also several new asymptotic and finite
upper bounds on the length function.

It was natural to analyze and survey the previously known
results, as well as presenting the new ones. In particular,
this was done for covering radius $R=2,3$. A survey of the most
used $q^{m}$-concatenating constructions is also given. It
should be noted that no surveys of nonbinary linear covering
codes have been recently published.

We also point out that new upper bounds on the length function are also new
upper bounds on the smallest possible sizes of saturating sets. More
generally, the new results and methods concerning small saturating sets in
projective spaces over finite fields that have been given in this paper,
such as the new concept of multifold strong blocking sets, seem to be of
independent interest.

\section*{Acknowledgments}

The authors thank Professor P. R. J. \"{O}sterg\aa rd for useful discussions
of problems on codes with covering radius three.


\begin{thebibliography}{99}
\bibitem{M-WSl} F. J.\ MacWilliams and N. J. A.\ Sloane, \emph{The theory of
error-correcting codes}. Amsterdam, The Netherlands: North-Holland, 1977.

\bibitem{Handbook-codes} V.\ S.\ Pless, W.\ C.\ Huffman, and R.~A.~Brualdi,
\textquotedblleft An introduction to algebraic codes,\textquotedblright\ in
\emph{Handbook of Coding Theory,} vol.\ 1, V.\ S.\ Pless, W.\ C.\ Huffman,
R.~A.~Brualdi, Eds.\ Amsterdam, The Netherlands: Elsevier, 1998, pp.\ 3-139.

\bibitem{Coh} G.\ Cohen, I.\ Honkala, S.\ Litsyn, and A.\
    Lobstein, \emph{ Covering Codes.} Amsterdam, The
    Netherlands: North-Holland, 1997.

\bibitem{Handbook-coverings} R.\ A.\ Brualdi, S.\ Litsyn, and V.\ S.\ Pless,
\textquotedblleft Covering Radius,\textquotedblright\ in \emph{Handbook of
Coding Theory,} vol.\ 1, V.\ S.\ Pless, W.\ C.\ Huffman, R.~A.~Brualdi,
Eds.\ Amsterdam, The Netherlands: Elsevier, 1998, pp.\ 755-826.

\bibitem{CohSurvey1985} G.\ Cohen, M. G. Karpovsky, H. F.
    Mattson, Jr., and J. R. Shatz,\textquotedblleft Covering
    radius - Survey and recent results,\textquotedblright\
    \emph{IEEE Trans. Inf.\ Theory,} vol. 31,
    no. 3, pp. 328-343, May 1985.

\bibitem{GrSl} R. L. Graham and N. J. A. Sloane,
    \textquotedblleft On the covering radius of
    codes,\textquotedblright\ \emph{IEEE Trans. Inf.\ Theory,}
     vol. 31, no. 3, pp. 385-401, May 1985.

\bibitem{KabPan} G. A. Kabatyansky and V. I. Panchenko,
    \textquotedblleft Unit sphere packings and coverings of the
    Hamming space,\textquotedblright\ \emph{Probl.
    Inf.\ Transm.}, vol. 24, no. 4, pp. 261-272, 1988.

\bibitem{BrPlWi} R.\ A.\ Brualdi, V.\ S.\ Pless, and R.\ M.\
    Wilson, \textquotedblleft Short codes with a given covering
    radius,\textquotedblright\ \emph{IEEE Trans.\ Inf.\
    Theory}, vol.\ 35, pp.\ 99-109, Jan.\ 1989.

\bibitem{Jan} H. Janwa, \textquotedblleft Some optimal codes from algebraic
geometry and their covering radii,\textquotedblright\ \emph{Europ.\ J.\
Combin.}, vol.\ 11, pp.\ 249-266, 1990.

\bibitem{DavPPI} A. A.\ Davydov, \textquotedblleft Construction of linear
covering codes,\textquotedblright\ \emph{Probl. Inf.\ Transm.}, vol. 26,
no.~4, pp.\ 317-331, 1990.

\bibitem{OstUpIEEE} P. R. J. \"{O}sterg\aa rd,
    \textquotedblleft Upper bounds for $q$-ary covering
    codes,\textquotedblright\ \emph{IEEE Trans.\ Inf.\ Theory},
     vol. 37, no. 3, pp. 660-664,
    May 1991; correction vol. 37, no. 6, p. 1738, Nov. 1991.

\bibitem{HonSubN91} I.\ S.\ Honkala, \textquotedblleft On
    $(k,t)$-subnormal covering codes,\textquotedblright\
    \emph{IEEE Trans.\ Inf.\ Theory}, vol. 37, no. 4,
    pp. 1203-1206, Jul. 1991.

\bibitem{Hon-length} ------, \textquotedblleft On lengthening of covering
codes,\textquotedblright\ \emph{Discrete Math.}, vol.\ 106--107, pp.\
291-295, 1992.

\bibitem{OstSubn92} P. R. J. \"{O}sterg\aa rd,
    \textquotedblleft Further results on $(k,t)$-subnormal
    covering codes,\textquotedblright\ \emph{IEEE Trans.\
    Inf.\ Theory}, vol. 38, no. 1, pp. 206-210, Jan.
    1992.

\bibitem{DavParis} A. A.\ Davydov, \textquotedblleft Construction of codes
with covering radius 2,\textquotedblright\ \emph{Algebraic Coding}
(G.~Cohen, S.\ Litsyn, A.\ Lobstein, G.\ Zemor, eds.), \emph{Lecture Notes
in Computer Science,} Springer-Verlag, New-York, vol. 573, pp.\ 23-31, 1992.

\bibitem{Rene-PhD} R.\ Struik, \emph{Covering codes}.\ Ph.D dissertation,
Eindhoven University of Technology, 1994.

\bibitem{DD-L} A.\ A.\ Davydov and A. Yu.
    Drozhzhina-Labinskaya, \textquotedblleft Constructions,
    families and tables of binary linear covering
    codes,\textquotedblright\ \emph{IEEE Trans.\ Inf.\ Theory,}
    vol. 40, no. 4, pp.\ 1270-1279, 1994.

\bibitem{Dav95} A. A.\ Davydov, \textquotedblleft Constructions
    and families of covering codes and saturated sets of points
    in projective geometry,\textquotedblright\ \emph{IEEE
    Trans.\ Inf.\ Theory}, vol.\ 41, no. 6,
    pp.\ 2071-2080, Nov.\ 1995.

\bibitem{davShumen} ------, \textquotedblleft On nonbinary
    linear codes with covering radius two,\textquotedblright\
    in \emph{Proc. 5th Int. Workshop Algebraic
    Combin. Coding Theory, }ACCT-V, Unicorn, Shumen,
    Bulgaria 1996, pp.~105-110.

\bibitem{Dav-Nonlin} ------, ``Constructions of
    nonlinear covering codes,'' \emph{IEEE
    Trans. Inf.\ Theory}, vol. 43, no. 5, pp.\
    1639-1647, Sep. 1997.

\bibitem{coc} J. C. Cock and P. R. J. \"{O}sterg\aa rd, \textquotedblleft
Ternary covering codes derived from BCH codes,\textquotedblright\ \emph{J.
Combin.\ Theory.\ Ser.\ A\/}, vol.\ 80, pp.\ 283-289, 1997.

\bibitem{BV} T. S. Baicheva and E. D. Velikova,
    \textquotedblleft Covering radii of ternary linear codes of
    small dimensions and codimensions,\textquotedblright\
    \emph{IEEE Trans.\ Inf.\ Theory}, vol. 43, no. 6,
    pp. 2057-2061, Nov. 1997; correction vol. 44, no. 5, p.
    2032, Sep. 1998.

\bibitem{ostArs} P. R. J. \"{O}sterg\aa rd, \textquotedblleft New
constructions for $q$-ary covering codes,\textquotedblright\ \emph{Ars
Combin.}, vol. 52, pp. 51-63, 1999.

\bibitem{DavNBCR2} A. A.\ Davydov, \textquotedblleft
    Constructions and families of nonbinary linear codes with
    covering radius~2,\textquotedblright\ \emph{IEEE Trans.\
    Inf.\ Theory} , vol.\ 45, no. 5, pp.\
    1679-1686, Jul. 1999.

\bibitem{JanwaMatt} H. Janwa and H. F. Mattson, Jr., \textquotedblleft Some
upper bounds on the covering radii of linear codes over $F_{q}$ and their
applications,\textquotedblright\ \emph{Des. Codes Cryptogr}, vol. 18, pp.
163-181, 1999.

\bibitem{Gab-NewtRad} E. M. Gabidulin and T. Kl\o ve, \textquotedblleft On
the Newton and covering radii of linear codes,\textquotedblright\ \emph{IEEE
Trans. Inf. Theory}, vol. 45, no. 7, pp. 2534-2536, Nov. 1999.

\bibitem{DavO1} A.\ A.\ Davydov and P.\ R.\ J. \"{O}sterg\aa rd,
\textquotedblleft New linear codes with covering radius 2 and odd
basis,\textquotedblright\ \emph{Des. Codes Cryptogr.}, vol. 16, pp. 29-39,
1999.

\bibitem{DavO3} ------, \textquotedblleft New quaternary linear codes with
covering radius 2,\textquotedblright\ \emph{Finite Fields Appl.}, vol. 6,
pp. 164-174, 2000.

\bibitem{DavO5} ------, \textquotedblleft Linear codes with
    covering radius $R=2,3$ and codimension $tR$,\textquotedblright\ \emph{IEEE
Trans.\ Inf.\ Theory}, vol. 47, no. 1, pp.~416-421, Jan. 2001.

\bibitem{Dav-2001-newconst} A. A.\ Davydov, \textquotedblleft New
constructions of covering codes,\textquotedblright\ \emph{Des. Codes
Cryptogr.}, vol. 22, pp.\ 305-316, 2001.

\bibitem{Baich2001} T. S. Baicheva, \textquotedblleft On the
    covering radius of ternary negacyclic codes with length up
    to 26,\textquotedblright\ \emph{ IEEE Trans. Inf. Theory},
    vol. 47, no. 1, pp. 413--416, Jan. 2001.

\bibitem{AshihBarg} A. Ashikhmin and A. Barg, \textquotedblleft Bounds on
the covering radius of linear codes,\textquotedblright\ \emph{Des. Codes
Cryptogr.}, vol. 27, pp. 261--269, 2002.

\bibitem{KaikRosTurku} M. K. Kaikkonen and P. Rosendahl,
    \textquotedblleft New covering codes from an ADS-like
    construction,\textquotedblright\ \emph{ IEEE Trans. Inf.
    Theory}, vol. 49, no. 7, pp. 1809-1812, Jul. 2003.

\bibitem{KrivSudVu} M. Krivelevich, B. Sudakov and V. H. Vu,
\textquotedblleft Covering codes with improved density,\textquotedblright\
\emph{IEEE Trans. Inf. Theory}, vol. 49, no. 7, pp. 1812-1815, Jul. 2003.

\bibitem{DMPIEEE} A.\ A.\ Davydov, S.\ Marcugini, and F.\
    Pambianco, \textquotedblleft Linear codes with covering
    radius 2,3 and saturating sets in projective
    geometry,\textquotedblright\ \emph{IEEE Trans.\
    Inf.\ Theory}, vol.\ 50, no. 3, pp.\
    537-541, Mar. 2004.

\bibitem{EtzMou} T. Etzion and B. Mounits, \textquotedblleft
    Quasi-perfect codes with small distance\textquotedblright ,
    \emph{IEEE Trans. Inf.\ Theory}, vol. 51, no. 11,
    pp. 3938-3946, Nov. 2005.

\bibitem{DFMP-IEEE-LO} A.\ A.\ Davydov, G. Faina, S.\
    Marcugini, and F.\ Pambianco, \textquotedblleft Locally
    optimal (nonshortening) linear covering codes and minimal
    saturating sets in projective spaces,\textquotedblright\
    \emph{IEEE Trans.\ Inf.\ Theory} , vol. 51,
    no. 12, pp. 4378-4387, Dec. 2005.

\bibitem{GiulPast} M. Giulietti and F. Pasticci,
    \textquotedblleft Quasi-perfect linear codes with minimum
    distance 4,\textquotedblright\ \emph{ IEEE Trans.
    Inf.\ Theory}, vol. 53, no. 5, pp.
    1928-1935, May 2007.

\bibitem{CohVar2007} G.\ Cohen and A. Vardy, \textquotedblleft Duality
between packings and coverings of the Hamming space,\textquotedblright\
\emph{Advances in Mathematics of Communications,} vol. 1, pp. 93-97, 2007.

\bibitem{BBDF-IEEE-2008} T. Baicheva, I. Bouyukliev, S.
    Dodunekov, and V. Fack, \textquotedblleft Binary and
    ternary linear quasi-perfect codes with small
    dimensions,\textquotedblright\ \emph{IEEE Trans.\
    Inf. Theory} , vol. 54, no. 9, pp. 4335-4339, Sep. 2008.

\bibitem{DavOprep} A.\ A.\ Davydov and P.\ R.\ J. \"{O}sterg\aa rd,
\textquotedblleft Linear codes with covering radius 3,\textquotedblright\ in
preparation.

\bibitem{Lobstein} A. Lobstein, \textquotedblleft Covering
    radius,\textquotedblright\ a bibliography. [Online].
    Available: http://www.infres.
    enst.fr/\symbol{126}lobstein/bib-a-jour.pdf.

\bibitem{GalKaba} F. Galand and G. Kabatiansky, \textquotedblleft
Information hiding by coverings,\textquotedblright\ in \emph{Proc. IEEE Inf.
Theory Workshop, }pp. 151--154, Paris, 2003.

\bibitem{BierbStegan} J. Bierbrauer and J. Fridrich,
    \textquotedblleft Constructing good covering codes for
    applications in Steganography,\textquotedblright\ in
    \emph{Lecture Notes in Computer Science, Trans.
    Data Hiding Multimedia Security III,} Springer-Verlag,
    vol. 4920, Y. Q. Shi, Ed., pp. 1-22, 2008.

\bibitem{Bo-Sz-Ti} E. Boros, T. Sz\H{o}nyi, and K. Tichler,
\textquotedblleft On defining sets for projective planes,\textquotedblright\
\emph{Discr. Math.,} vol. 303, pp. 17-31, 2005.

\bibitem{KKKPS} G. Kiss, I. Kov\'{a}cs, K. Kutnar, J. Ruff and P. \v{S}parl,
\textquotedblleft A note on a geometric construction of large Cayley graphs
of given degree and diameter,\textquotedblright\ submitted.

\bibitem{PartSumQuer} C. T. Ho, J. Bruck, and R. Agrawal, \textquotedblleft
Partial-sum queries in OLAP data cubes using covering
codes,\textquotedblright\ \emph{IEEE Trans. Computers, }vol. 47, no. 12, pp.
1326-1340, Dec. 1998.

\bibitem{Hats-on-line} S. Aravamuthan and S. Lodha, \textquotedblleft
Covering codes for hats-on-a-line,\textquotedblright\ \emph{Electronic J.
Combin}., vol. 13, \#R21, 2006.

\bibitem{Identif} G. Exoo, V. Junnila, T. Laihonen, S. Ranto,
    \textquotedblleft Constructions for identifying
    codes,\textquotedblright\ in \emph{Proc. XI Int. Workshop
    Algebraic Comb. Coding Theory}, ACCT2008, Pamporovo,
    Bulgaria, Jun. 2008, pp.
    92-98.  [Online]. Available: http://www.moi.math.bas.bg/acct2008/b16.pdf

\bibitem{Hirs} J.\ W.\ P.\ Hirschfeld, \emph{Projective geometries over
finite fields,} second edition, Oxford University Press, Oxford, 1998.

\bibitem{Hirs1} ------, \textquotedblleft Maximum sets in finite projective
spaces\textquotedblright , in \emph{Surveys in Combinatorics}, E. K. Lloyd,
Ed., London Math.\ Soc.\ Lecture Note Series, vol.\ 82, Cambridge, UK:
Cambridge Univ.\ Press, pp.\ 55-76, 1983.

\bibitem{HS1} J. W. P. Hirschfeld and L. Storme, \textquotedblleft The
packing problem in statistics, coding theory and finite projective
spaces,\textquotedblright\ \emph{J. Statist. Planning Infer.}, vol. 72, pp.
355-380, 1998.

\bibitem{HS} ------, \textquotedblleft The packing problem in statistics,
coding theory and finite projective spaces: update 2001,\textquotedblright\
Blokhuis, A. (ed.) et al., \emph{Finite geometries}. Proceedings of the
fourth Isle of Thorns conference, Brighton, UK, April 2000. Dordrecht:
Kluwer Academic Publishers. Dev. Math., vol. 3, pp. 201-246, 2001.

\bibitem{Sz-survey89} T.~Sz{\H{o}}nyi, \textquotedblleft
    Complete arcs in Galois planes: a
    survey,\textquotedblright\ \emph{Quaderni del seminario di
    Geometrie Combinatorie,} vol. 94, Dipartimento di
    Matematica G. Castelnuovo, Univ. di Roma La Sapienza, 1989.

\bibitem{Lang2000DM} I. N Landjev, \textquotedblleft Linear
    codes over finite fields and finite projective
    geometries,\textquotedblright\ \emph{ Discr. Math.}, vol.
    213, pp. 211-244, 2000.

\bibitem{Bartocci} U. Bartocci, \textquotedblleft $k$-insiemi densi in piani
di Galois,\textquotedblright\ \emph{Boll. Un. Mat. Ital.} D vol. 2, pp.
71-77, 1983.

\bibitem{Ughi} E.\ Ughi, \textquotedblleft Saturated configurations of
points in projective Galois spaces,\textquotedblright\ \emph{Europ J.\
Combin.,} vol.\ 8, pp.\ 325-334, 1987.

\bibitem{DavO4} A.\ A.\ Davydov and P.\ R.\ J. \"{O}sterg\aa rd,
\textquotedblleft On saturating sets in small projective
geometries,\textquotedblright\ \emph{European J. Combin.}, vol. 21, pp.
563-570, 2000.

\bibitem{MPAustr} S.\ Marcugini and F.\ Pambianco,
    \textquotedblleft Minimal 1-saturating sets in $PG(2,q),$
    $q\leq 16,$\textquotedblright\ \emph{ Austral. J.
    Combin.}, vol. 28, pp. 161-169, 2003.

\bibitem{DMP-JCTA} A.\ A.\ Davydov, S.\ Marcugini, and F.\ Pambianco,
\textquotedblleft On saturating sets in projective
spaces,\textquotedblright\ \emph{J.\ Combin.\ Theory, Ser.\ A},\emph{\ }
vol.\ 103, pp.\ 1-15, 2003.

\bibitem{kovS92} S.~J. Kov\'{a}cs,\textquotedblleft Small saturated sets in
finite projective planes,\textquotedblright\ \emph{Rend. Mat.}, Ser. VII,
vol. 12, pp. 157--164, 1992.

\bibitem{GiuTor-Ars04} M. Giulietti and F. Torres, \textquotedblleft On
dense sets related to plane algebraic curves,\textquotedblright\ \emph{Ars.
Combin}., vol. 72, pp. 33-40, 2004.

\bibitem{Giul-plane} M. Giulietti, \textquotedblleft On small dense sets in
Galois planes,\textquotedblright\ \emph{Electronic J. Combin.}, vol. 14,
Research Paper 75, 13 pp., 2007.

\bibitem{DGMP-ACCT2008} A. A. Davydov, M. Giulietti, S.
    Marcugini and F. Pambianco, \textquotedblleft Linear
    covering codes over nonbinary finite
    fields,\textquotedblright\ in \emph{Proc. XI Int. Workshop
    Algebraic Comb. Coding Theory}, ACCT2008, Pamporovo,
    Bulgaria, Jun. 2008, pp. 70-75. [Online]. Available:
    http://www.moi.math.bas.bg/acct2008/b12.pdf

\bibitem{DGMP-Petersb2008} A. A. Davydov, M. Giulietti, S.
    Marcugini and F. Pambianco, \textquotedblleft Linear
    Covering Codes of Radius 2 and 3,\textquotedblright\ in
    \emph{Proc. Workshop \textquotedblleft Coding Theory Days
    in St. Petersburg\textquotedblright }, St. Petersburg,
    Russia, Oct. 2008, pp. 12-17. [Online]. Available:
    http://k36.org/codingdays/proceedings.pdf

\bibitem{balS96} S.~Ball and A.~Blokhuis, \textquotedblleft On the size of a
double blocking set in {$PG(2,q)$},\textquotedblright\ \emph{Finite Fields
Appl.}, vol. 2, pp. 125-137, 1996.

\bibitem{BallHirs} S.~Ball and J.\ W.\ P.\ Hirschfeld, \textquotedblleft
Bounds on $(n,r)$-arcs and their application to linear
codes,\textquotedblright\ \emph{Finite Fields Appl.,} vol. 11, pp. 326-336,
2005.

\bibitem{DFMP-JG} A.\ A.\ Davydov, G. Faina, S.\ Marcugini, and F.\
Pambianco, \textquotedblleft Computer search in projective planes for the
sizes of complete arcs,\textquotedblright\ \emph{J.\ Geom.}, vol. 82, pp.
50-62, 2005.

\bibitem{Gi2} M. Giulietti, \textquotedblleft Small complete
    caps in $ PG(N,q) $, $q$ even,\textquotedblright\ \emph{J.
    Comb. Des.}, vol. 15 , pp. 420--436, 2007.

\bibitem{DGMP-Submit} A. A.\ Davydov, M.\ Giulietti, S.\ Marcugini and F.\
Pambianco, \textquotedblleft New inductive constructions of complete caps in
$PG(N,q)$, $q$ even,\textquotedblright\ submitted.

\bibitem{Ost2002} P.\ R.\ J. \"{O}sterg\aa rd, \textquotedblleft Classifying
Subspaces of Hamming Spaces," \emph{Des.\ Codes Crypt.}, vol.\ 27, pp.\
297-305, 2002.
\end{thebibliography}
\end{document}